\documentclass[10pt,twoside]{article}
\usepackage{mathrsfs}
\usepackage{amsmath}
\usepackage{amssymb}
\usepackage{fancyhdr}
\usepackage{latexsym}
\usepackage{bbding}
\usepackage{mathrsfs}
\usepackage{wasysym}
\usepackage{multicol,graphics}
\usepackage{cite}

\newtheorem{theorem}{Theorem}[section]
\newtheorem{lemma}[theorem]{Lemma}

\newtheorem{definition}[theorem]{Definition}

\numberwithin{equation}{section}

\allowdisplaybreaks

\makeatletter
\begin{document}
\title{{Global Existence of Large Solutions for the 3D incompressible Navier--Stokes--Poisson--Nernst--Planck Equations}\footnote{Corresponding author: jihzhao@163.com (J. Zhao).}}
\author{Jihong Zhao, Ying Li\\
[0.2cm] {\small School of Mathematics and Information Science, Baoji University of Arts and Sciences,}\\
[0.2cm] {\small  Baoji, Shaanxi 721013,  China}\\
[0.2cm] {\small E-mail: jihzhao@163.com, yingl723@163.com}}
\date{}
\maketitle

\begin{abstract}
This work is concerned with the global existence of large solutions to the three-dimensional dissipative fluid-dynamical
model, which is a strongly coupled nonlinear nonlocal system characterized by the incompressible
Navier--Stokes--Poisson--Nernst--Planck equations. Making full use of the algebraic structure of the system, we obtain the global existence of solutions without smallness assumptions imposed on the third component of the initial velocity field and the summation of initial densities of  charged species. More precisely, we prove that there
exist two positive constants $c_{0}, C_{0}$ such that if the initial data satisfies
\begin{align*}
  \big(\|u_{0}^{h}\|_{\dot{B}^{-1+\frac{3}{p}}_{p,1}}+\|N_{0}-P_{0}\|_{\dot{B}^{-2+\frac{3}{q}}_{q,1}}
  \big)
  &\exp\Big\{C_{0}\big(\|u_{0}^{3}\|_{\dot{B}^{-1+\frac{3}{p}}_{p,1}}^{2}+(\|N_{0}+P_{0}\|_{\dot{B}^{-2+\frac{3}{r}}_{r,1}}+1)\\
   &\times\exp\big\{C_{0}\|u_{0}^{3}\|_{\dot{B}^{-1+\frac{3}{p}}_{p,1}}\big\}+1\big)\Big\}
   \leq
  c_{0},
\end{align*}
then the incompressible Navier--Stokes--Poisson--Nernst--Planck equations admits a unique global solution.


\end{abstract}

\section{Introduction}
In this paper, we study the Cauchy problem of three-dimensional (3D) incompressible Navier--Stokes--Poisson--Nernst--Planck equations:
\begin{equation}\label{eq1.1}
\begin{cases}
  \partial_{t} u+u\cdot\nabla u-\Delta
  u+\nabla \pi=\Delta
  \phi\nabla\phi,\\
  \nabla\cdot u=0,\\
  \partial_{t} N+u\cdot \nabla
  N=\nabla\cdot(\nabla N-N\nabla \phi),\\
  \partial_{t} P+u\cdot \nabla
  P=\nabla\cdot(\nabla P+P\nabla \phi),\\
 \Delta \phi=N-P,\\
 (u, N, P)|_{t=0}=(u_0, N_0, P_0),
\end{cases}
\end{equation}
where $(x,t)\in \mathbb{R}^3\times\mathbb{R}_{+}$, $u=(u^{1}(x,t),u^{2}(x,t),u^{3}(x,t))$ and $\pi=\pi(x,t)$ stand for the
velocity field and the pressure of the incompressible fluid, respectively,  $N=N(x,t)$ and $P=P(x,t)$ stand for the densities of a negatively and positively charged species,
respectively, and $\phi=\phi(x,t)$ is the electrostatic potential caused by the charged species. For the sake of simplicity of presentation, we have assumed that the fluid density, viscosity, charge mobility and dielectric constant are unity.

The first two equations of \eqref{eq1.1} are the momentum
conservation and the mass conservation equations of the
incompressible flow, and the right-hand side
term in the momentum equations is the Lorentz force, which exhibits
$\Delta\phi\nabla\phi=\nabla\cdot\sigma$,
and the electric stress $\sigma$ is a rank one tensor plus a
pressure, for $i,j=1,2,3$,
\begin{equation}\label{eq1.2}
  [\sigma]_{ij}=\big(\nabla\phi\otimes\nabla\phi-\frac{1}{2}|\nabla\phi|^{2}I\big)_{ij}
  =\partial_{x_{i}}\phi\partial_{x_{j}}\phi-\frac{1}{2}|\nabla\phi|^{2}\delta_{ij},
\end{equation}
where $I$ is $3\times3$ identity matrix,  $\delta_{ij}$ is the
Kronecker symbol, and $\otimes$ denotes the tensor product. The
electric stress $\sigma$ stems from the balance of kinetic energy
with electrostatic energy via the least action principle (cf.
\cite{RLZ07}). The third and fourth equations of \eqref{eq1.1} model the balance
between diffusion and convective transport of the charged species by the
flow and the electric fields, and the fifth equation of \eqref{eq1.1} is the Poisson equation for the
electrostatic potential $\phi$, where the right-hand side is the net charge density.

The system \eqref{eq1.1} was first introduced by Rubinstein
\cite{R90}, which is capable of describing electro-chemical and
fluid-mechanical transport throughout the cellular environment. At
the present time, modeling of electro-diffusion in electrolytes is a
problem of major scientific interest, it finds that such model has a
wide applications in biology (ion channels), chemistry
(electro-osmosis) and pharmacology (transdermal iontophoresis), we
refer the readers to see \cite{BW04, EN00, ES05, SK04} for more detailed  applications of the system
\eqref{eq1.1} in electro-hydrodynamics, and \cite{LCJS081, LCJS082, LMCJS08} for the computational simulations.

The mathematical analysis of the system \eqref{eq1.1} was initiated by  Jerome \cite{J02}, where  the author established a local existence--uniqueness theory of the system \eqref{eq1.1}
 in the Kato's analytical  semigroup framework. Since the right-hand side term $\Delta\phi\nabla\phi$ has a nice algebraic structure \eqref{eq1.2},  it can be regarded that $\nabla\phi$ plays the same role as the velocity field $u$. Based on this observation,  by using the Hardy--Littlewood--Sobolev inequality and the Sobolev embeddings $W^{1,\frac{3}{2}}(\mathbb{R}^{3})\hookrightarrow L^{3}(\mathbb{R}^{3})$ and $\dot{B}^{-1+\frac{3}{q}}_{q,1}(\mathbb{R}^{3})\hookrightarrow\dot{B}^{-1+\frac{3}{p}}_{p,1}(\mathbb{R}^{3})$
with $1\leq q\le p\leq \infty$,  Zhao--Deng--Cui \cite{ZDC10, ZDC11} established local well-posedness and global
well-posedness with small initial
data of the system \eqref{eq1.1} in critical Lebesgue spaces and Besov spaces under the heat semigroup framework. For more analytical results concerning about the global existence of (large) weak solutions and (small) mild solutions, convergence rate estimates to
stationary solutions of time-dependent solutions and other related topics we refer the readers to see \cite{CI19, FS17, J11, JS09, R09, S09, WJL19, WJ21, ZL17} and references therein.

In order to give a better explanation of our main results, we let $v:=N-P$, $w:=N+P$, and the system
\eqref{eq1.1} turns into
\begin{equation}\label{eq1.3}
\begin{cases}
  \partial_{t} u+u\cdot\nabla u-\Delta
  u+\nabla \pi=-v\nabla(-\Delta)^{-1}v,\\
  \nabla\cdot u=0,\\
  \partial_{t} v+u\cdot \nabla
  v=\nabla\cdot(\nabla v+w\nabla(-\Delta)^{-1}v),\\
  \partial_{t} w+u\cdot \nabla
  w=\nabla\cdot(\nabla w+v\nabla(-\Delta)^{-1}v),\\
  (u, v, w)|_{t=0}=(u_0, v_0, w_0),
\end{cases}
\end{equation}
where $v_{0}=N_{0}-P_{0}$ and $w_{0}=N_{0}+P_{0}$. Moreover, let $u:=(u^{h}, u^{3})$, where $u^{h}:=(u^{1},u^{2})$ and $u^{3}$ denote the horizontal components and vertical component of the velocity field $u$, respectively. Using the divergence free condition $\nabla\cdot{u}=0$, it is easy to verify that the vertical component $u^{3}$ satisfies the following equation:
\begin{equation}\label{eq1.4}
  \partial_{t}u^{3}-\Delta
  u^{3}= -\operatorname{div}_{h}(u^{h}u^{3})+2u^{3}\operatorname{div}_{h}u^{h}
  -\partial_{3}\pi-v\partial_{3}(-\Delta)^{-1}v,
\end{equation}
which reveals that the equation on the vertical component $u^{3}$ is actually a linear equation with coefficients depending on the
horizontal components $u^{h}$ and the density function $v$.  Based on this observation, the authors in \cite{ZZL15} proved that under the conditions $1<q\leq p<6$ and
$\frac{1}{p}+\frac{1}{q}\geq\frac{2}{3}$, there exist two
positive constants $c_{0}$ and $C_{0}$  such that if
the initial data $(u_{0},v_{0}, w_{0})$ satisfies
\begin{equation}\label{eq1.5}
  \big(\|u_{0}^{h}\|_{\dot{B}^{-1+\frac{3}{p}}_{p,1}}+\|v_{0}\|_{\dot{B}^{-2+\frac{3}{q}}_{q,1}}
  +\|w_{0}\|_{\dot{B}^{-2+\frac{3}{q}}_{q,1}}\big)
  \exp\Big\{C_{0}(\|u_{0}^{3}\|_{\dot{B}^{-1+\frac{3}{p}}_{p,1}}^{2}+1)\Big\}\leq
  c_{0},
\end{equation}
then system \eqref{eq1.3} admits a unique global
solution. On the other hand,  observing that the fourth equation of \eqref{eq1.3} is a linear equation for $w$ with coefficients depending on the velocity field $u$ and the density function $v$, which may suggest us that we do not need to impose any smallness condition on initial data $w_{0}$ to ensure global existence of solutions. Indeed, Ma \cite{M18} showed that under the conditions $1\leq  p<\infty$, $1\leq q<6$, $q\leq 2p$ and
$\frac{1}{p}-\frac{1}{q}<\frac{1}{3}<\frac{1}{p}+\frac{1}{q}$, there exist two
positive constants $c_{0}$ and $C_{0}$ such that if
the initial data $(u_{0},v_{0}, w_{0})$ satisfies
\begin{equation}\label{eq1.6}
  \big(\|u_{0}\|_{\dot{B}^{-1+\frac{3}{p}}_{p,1}}+\|v_{0}\|_{\dot{B}^{-2+\frac{3}{q}}_{q,1}}
  \big)
  \exp\Big\{C_{0}\|w_{0}\|_{\dot{B}^{-2+\frac{3}{q}}_{q,1}}\Big\}\leq
  c_{0},
\end{equation}
then system \eqref{eq1.3} still has a unique global solution.

Motivated by the above global existence results in \cite{M18, ZZL15}, in this paper, we
aim at relaxing the smallness conditions imposed on the vertical component of initial velocity field and the summation of the initial densities of charged species
such that system \eqref{eq1.3}  still has a unique global solution. Considering the algebraic structures of the nonlinear terms in \eqref{eq1.3}, by \cite{I11}, the nonlinear term $v\nabla(-\Delta)^{-1}v$ has a nice symmetric structure as
\begin{align}\label{eq1.7}
\partial_{x_{i}}v\partial_{x_{i}}(-\Delta)^{-1}v=&\frac{1}{2}\partial_{x_{i}}(-\Delta) \{((-\Delta)^{-1}v)(\partial_{x_{i}}(-\Delta)^{-1}v)\}\nonumber\\
&+\partial_{x_{i}}\nabla\cdot \{((-\Delta)^{-1}v)(\nabla\partial_{x_{i}}(-\Delta)^{-1}v)\}\nonumber\\
&+\frac{1}{2}\partial_{x_{i}}^{2} \{((-\Delta)^{-1}v)v\}.
\end{align}
This enables us to treat the
equation of $w$ in a weaker Besov space $\dot{B}^{-2+\frac{3}{r}}_{r,1}$ with $1\leq r<\infty$. However, the nonlinear coupled term $w\nabla (-\Delta)^{-1}v$ has lack
of such a symmetric structure, which can not  exhibit such a good expression
as \eqref{eq1.7} and prevents us to obtain good estimates for the equation $v$ in a weaker Besov spaces.   These observations essentially indicate that the
difference of charged densities $v$ plays more important role than the summation of charged densities $w$ in mathematical analysis of the system \eqref{eq1.3}.   Based on these careful observations, by using analytical methods in \cite{M18, ZZL15}, we intend to prove the global existence of solutions under
the assumptions that the horizontal components of the velocity field and the difference of initial densities are small while the vertical component of the velocity field and the total initial densities could be chosen suitable large. Moreover, we consider the functional space of solutions of the system \eqref{eq1.3} with initial data $v_{0}$ and $w_{0}$ belonging to the different low regularity Besov spaces with different regularity and integral indices, which can indicate more specific coupling relations between the difference and the summation of  negatively and positively charged densities.

Now we state our main results as  follows.

\begin{theorem}\label{th1.1}  Let $p,q,r$ be three positive numbers such that $1\leq p, q, r<\infty$, and satisfy the following conditions:
\begin{equation*}
  \frac{1}{q}-\frac{1}{p}\geq -\min\{\frac{1}{3}, \frac{1}{2p}\},\ \ \
  \max\{\frac{1}{q}-\frac{1}{r}, \frac{1}{r}-\frac{1}{q}\}<\frac{1}{3}<\min\{\frac{1}{p}+\frac{1}{q}, \frac{1}{p}+\frac{1}{r}, \frac{1}{q}+\frac{1}{r}\}.
\end{equation*}
Then for any $u_{0}\in
\dot{B}^{-1+\frac{3}{p}}_{p,1}(\mathbb{R}^{3})$ with $\nabla\cdot
u_{0}=0$, $v_{0}\in
\dot{B}^{-2+\frac{3}{q}}_{q,1}(\mathbb{R}^{3})$ and $ w_{0}\in
\dot{B}^{-2+\frac{3}{r}}_{r,1}(\mathbb{R}^{3})$, there exists
$T>0$ such that the system \eqref{eq1.3} admits a unique solution $(u,v,w)$ on $[0,T]$ satisfying
\begin{equation}\label{eq1.8}
\begin{cases}
  u\in C([0,T],\dot{B}^{-1+\frac{3}{p}}_{p,1}(\mathbb{R}^{3}))\cap\mathcal{L}^{\infty}(0, T;
  \dot{B}^{-1+\frac{3}{p}}_{p,1}(\mathbb{R}^{3}))\cap  L^{1}(0, T;
  \dot{B}^{1+\frac{3}{p}}_{p,1}(\mathbb{R}^{3})),\\
  v\in C([0,T],\dot{B}^{-2+\frac{3}{q}}_{q,1}(\mathbb{R}^{3}))\cap\mathcal{L}^{\infty}(0, T;
  \dot{B}^{-2+\frac{3}{q}}_{q,1}(\mathbb{R}^{3}))\cap  L^{1}(0, T;
  \dot{B}^{\frac{3}{q}}_{q,1}(\mathbb{R}^{3})),\\
   w\in C([0,T],\dot{B}^{-2+\frac{3}{r}}_{r,1}(\mathbb{R}^{3}))\cap\mathcal{L}^{\infty}(0, T;
  \dot{B}^{-2+\frac{3}{r}}_{r,1}(\mathbb{R}^{3}))\cap  L^{1}(0, T;
  \dot{B}^{\frac{3}{r}}_{r,1}(\mathbb{R}^{3})).
\end{cases}
\end{equation}
Besides, there exists a positive constant $\varepsilon$ such that if the initial data satisfies
\begin{equation*}
\|(u_{0}, v_{0},
w_{0})\|_{\dot{B}^{-1+\frac{3}{p}}_{p,1}\times\dot{B}^{-2+\frac{3}{q}}_{q,1}\times\dot{B}^{-2+\frac{3}{r}}_{r,1}}
\leq \varepsilon,
\end{equation*}
then the above assertion holds for $T=\infty$, i.e., the solution $(u, v, w)$ is global.
\end{theorem}

\begin{theorem}\label{th1.2}
 Let $p,q,r$ be three positive numbers such that $1<p<6$, $1<q, r<\infty$, and satisfy the following conditions:
\begin{equation*}
  \frac{1}{q}-\frac{1}{p}\geq -\min\{\frac{1}{3}, \frac{1}{2p}\},\ \ \
  \max\{\frac{1}{q}-\frac{1}{r}, \frac{1}{r}-\frac{1}{q}\}<\frac{1}{3}<\min\{\frac{1}{p}+\frac{1}{q}, \frac{1}{p}+\frac{1}{r}, \frac{1}{q}+\frac{1}{r}\}.
\end{equation*}
Then for any $u_{0}=(u_{0}^{h},u_{0}^{3})\in
\dot{B}^{-1+\frac{3}{p}}_{p,1}(\mathbb{R}^{3})$ with $\nabla\cdot
u_{0}=0$, $v_{0}\in
\dot{B}^{-2+\frac{3}{q}}_{q,1}(\mathbb{R}^{3})$ and $ w_{0}\in
\dot{B}^{-2+\frac{3}{r}}_{r,1}(\mathbb{R}^{3})$, there exist two
positive constants $c_{0}$ and $C_{0}$ such that if
the initial data $(u_{0},v_{0}, w_{0})$ satisfies
\begin{equation}\label{eq1.9}
  \big(\|u_{0}^{h}\|_{\dot{B}^{-1+\frac{3}{p}}_{p,1}}+\|v_{0}\|_{\dot{B}^{-2+\frac{3}{q}}_{q,1}}
  \big)
  \exp\Big\{C_{0}\big(\|u_{0}^{3}\|_{\dot{B}^{-1+\frac{3}{p}}_{p,1}}^{2}+(\|w_{0}\|_{\dot{B}^{-2+\frac{3}{r}}_{r,1}}+1)
   \exp\big\{C_{0}\|u_{0}^{3}\|_{\dot{B}^{-1+\frac{3}{p}}_{p,1}}\big\}+1\big)\Big\}\leq
  c_{0},
\end{equation}
then the system \eqref{eq1.3} admits a unique global
solution $(u,v,w)$ satisfying \eqref{eq1.8}
\end{theorem}

\noindent\textbf{Remark 1.1}   The initial condition \eqref{eq1.9} exhibits that the initial data $u_{0}^{3}$ and $w_0$ can be taken large as long as we take the initial data $u_{0}^{h}$ and $v_0$ small enough compared with the size of $u_{0}^{3}$ and  $w_0$, which we can still get the global existence of solutions to the system \eqref{eq1.3}. Back to the original system \eqref{eq1.1},  the condition \eqref{eq1.9} is equivalent to the following condition:
\begin{equation}\label{eq1.10}
  \big(\|u_{0}^{h}\|_{\dot{B}^{-1+\frac{3}{p}}_{p,1}}+\|N_{0}-P_{0}\|_{\dot{B}^{-2+\frac{3}{q}}_{q,1}}
  \big)
  \exp\Big\{C_{0}\big(\|u_{0}^{3}\|_{\dot{B}^{-1+\frac{3}{p}}_{p,1}}^{2}+(\|N_{0}+P_{0}\|_{\dot{B}^{-2+\frac{3}{r}}_{r,1}}+1)
   \exp\big\{C_{0}\|u_{0}^{3}\|_{\dot{B}^{-1+\frac{3}{p}}_{p,1}}\big\}+1\big)\Big\}\leq
  c_{0},
\end{equation}
 thus Theorem \ref{th1.2} implies global existence of solutions for the system \eqref{eq1.1} with  only requiring the horizontal components of the initial velocity field and the difference of initial negatively and positively charged densities are small enough.

\noindent\textbf{Remark 1.2}  The specific coupled relation between $v$ and $w$ was indicated by the condition $\max\{\frac{1}{q}-\frac{1}{r}, \frac{1}{r}-\frac{1}{q}\}<\frac{1}{3}<\frac{1}{q}+\frac{1}{r}$, which tells us that the regularity of solution $v$ or $w$ can be taken beyond the regularity index $-\frac{3}{2}$, but one can not take both of them less than $-\frac{3}{2}$ at the same time.  Indeed, the regularity of $v$ can be taken much weaker as long as the regularity of $w$ is not that much weaker, i.e.,  $q$ can be taken large enough as long as we take $r$ closing to $3$ such that the condition $\max\{\frac{1}{q}-\frac{1}{r}, \frac{1}{r}-\frac{1}{q}\}<\frac{1}{3}<\frac{1}{q}+\frac{1}{r}$ holds.  Hence, Theorem \ref{th1.2} can be regarded as an extension of global existence results in \cite{ZDC10, ZZL15, M18}, where the global existence of solutions with small initial data was proved in critical Besov spaces with the same regularity and integral indices for $v$ and $w$, and the regularity index must less than $-\frac{3}{2}$.

This paper is organized as follows. In section 2, we first introduce definitions of the homogeneous Besov spaces and the Chemin--Lerner mixed time-space spaces based on the Littlewood--Paley dyadic decomposition theory, then review some known bilinear estimates which used frequently in the proofs of Theorems \ref{th1.1} and \ref{th1.2}.  In Section 3,  we first establish two crucial nonlinear estimates of the pressure $\pi$, then derive the desired estimates of $u^{h}$ and $u^{3}$ by using the weighted Chemin--Lerner type norm; while in Section 4, we derive the desired estimates of
$v$ and $w$. Finally in the last section, we complete the proofs of Theorems \ref{th1.1} and \ref{th1.2}.

\section{Preliminaries}

Throughout this paper,  $C$ and $C_{i}$ ($i=1,2, \cdots$) stand for the generic harmless constants. For brevity, we shall use the notation
$f\lesssim g$ instead of $f\le Cg$, and $f\approx g$ means that $f\lesssim g$
and $g\lesssim f$.  For any Banach spaces $\mathcal{X}$ and $\mathcal{Y}$, $f\in\mathcal{X}$ and $g\in\mathcal{Y}$, we write $\|(f,g)\|_{\mathcal{X}\times\mathcal{Y}}:=
\|f\|_{\mathcal{X}}+\|g\|_{\mathcal{Y}}$. For all $T>0$ and $\rho\in[1,\infty]$,  we denote by $C([0,T],\mathcal{X})$ the set of continuous functions on $[0,T]$ with
values in $\mathcal{X}$, and denote by $L^{\rho}(0,T;\mathcal{X})$ the set of
measurable functions on $[0,T]$ with values in $\mathcal{X}$ such that
$t\rightarrow \|f(t)\|_{\mathcal{X}}$ belongs to $L^{\rho}(0,T)$.

\subsection{Littlewood--Paley dyadic decomposition and Besov spaces}

Let us briefly recall the Littlewood--Paley dyadic decomposition theory and the stationary/time dependent Besov spaces for convenience.  More details may be found in Chap. 2 and Chap. 3 in the book \cite{BCD11}. Let $\mathcal{S}(\mathbb{R}^{3})$ be the Schwartz class of rapidly
decreasing functions, and $\mathcal{S}'(\mathbb{R}^{3})$ the space of tempered  distributions.  Choose a smooth radial non-increasing funciton $\chi$  with $\operatorname{Supp}\chi\subset B(0,\frac{4}{3})$ and $\chi\equiv1$ on $B(0,\frac{3}{4})$. Set $\varphi(\xi)=\chi(\frac{\xi}{2})-\chi(\xi)$. It is not difficult to check that
 $\varphi$ is supported in the shell $\{\xi\in\mathbb{R}^{3},\ \frac{3}{4}\leq
|\xi|\leq \frac{8}{3}\}$, and
\begin{align*}
   \sum_{j\in\mathbb{Z}}\varphi(2^{-j}\xi)=1 \ \ \text{for}\ \  \xi\in\mathbb{R}^{3}\backslash\{0\}.
\end{align*}
Let $h=\mathcal{F}^{-1}\varphi$. Then  for any $f\in\mathcal{S}'(\mathbb{R}^{3})$, the homogeneous dyadic blocks $\Delta_{j}$ ($j\in\mathbb {Z}$) are defined by
\begin{align*}
 \Delta_{j}f(x): =\varphi(2^{-j}D)f(x)=2^{3j}\int_{\mathbb{R}^{3}}h(2^{j}y)f(x-y)dy.
\end{align*}
Let $\mathcal{S}'_{h}(\mathbb{R}^{3})$ be the space of tempered distribution $f\in\mathcal{S}'(\mathbb{R}^{3})$ such that
$$
\lim_{j\rightarrow -\infty} S_{j}f(x)=0,
$$
where $S_{j}f$ ($j\in\mathbb {Z}$) stands for the low frequency cut-off defined by $ S_{j}f: =\chi(2^{-j}D)f$.
Then one has the unit decomposition for any tempered distribution $f\in\mathcal{S}'_{h}(\mathbb{R}^{3})$:
\begin{align}\label{eq2.1}
 f=\sum_{j\in\mathbb{Z}}\Delta_{j}f.
\end{align}
The above homogeneous dyadic block $\Delta_{j}$ and
the partial summation operator $S_{j}$  satisfy the following quasi-orthogonal properties: for any $f, g\in\mathcal{S}'(\mathbb{R}^{3})$, one has
\begin{align}\label{eq2.2}
  \Delta_{i}\Delta_{j}f\equiv0\ \ \ \text{if}\ \ \ |i-j|\geq 2\ \ \ \text{and}\ \ \
  \Delta_{i}(S_{j-1}f\Delta_{j}g)\equiv0 \ \ \ \text{if}\ \ \ |i-j|\geq 5.
\end{align}
Moreover, using Bony's homogeneous paraproduct decomposition (cf. \cite{B81}),  one can formally split the product of two temperate distributions $f$ and $g$ as follows:
\begin{equation}\label{eq2.3}
  fg=T_{f}g+T_{g}f+R(f,g),
\end{equation}
where the
paraproduct between $f$ and $g$ is defined by
\begin{equation*}
  T_{f}g:=\sum_{j\in\mathbb{Z}}S_{j-1}f\Delta_{j}g=\sum_{j\in\mathbb{Z}}\sum_{k\leq j-2}\Delta_{k}f\Delta_{j}g,
\end{equation*}
and the remaining term is defined by
\begin{equation*}
  R(f,g):=\sum_{j\in\mathbb{Z}}\Delta_{j}f\widetilde{\Delta_{j}}g \ \
  \text{and}\ \
  \widetilde{\Delta_{j}}:=\Delta_{j-1}+\Delta_{j}+\Delta_{j+1}.
\end{equation*}

Based on those dyadic blocks, the homogeneous Besov spaces can be defined as follows:
\begin{definition}\label{de2.1}
 For any $s\in \mathbb{R}$, $1\leq p,r\leq\infty$ and $f\in\mathcal{S}'(\mathbb{R}^{3})$, we set
\begin{equation*}
  \|f\|_{\dot{B}^{s}_{p,r}}:= \begin{cases} \left(\sum_{j\in\mathbb{Z}}2^{srj}\|\Delta_{j}f\|_{L^{p}}^{r}\right)^{\frac{1}{r}}
  \ \ &\text{for}\ \ 1\leq r<\infty,\\
  \sup_{j\in\mathbb{Z}}2^{sj}\|\Delta_{j}f\|_{L^{p}}\ \
  &\text{for}\ \
  r=\infty,
 \end{cases}
\end{equation*}
and the homogeneous Besov
space $\dot{B}^{s}_{p,r}(\mathbb{R}^{3})$ is defined by
\begin{itemize}
\item For $s<\frac{3}{p}$ (or $s=\frac{3}{p}$ if $r=1$), we define
\begin{equation*}
  \dot{B}^{s}_{p,r}(\mathbb{R}^{3}):=\Big\{f\in \mathcal{S}_{h}'(\mathbb{R}^{3}):\ \
  \|f\|_{\dot{B}^{s}_{p,r}}<\infty\Big\}.
\end{equation*}
\item If $k\in\mathbb{N}$ and $\frac{3}{p}+k\leq s<\frac{3}{p}+k+1$ (or $s=\frac{3}{p}+k+1$ if $r=1$), then $\dot{B}^{s}_{p,r}(\mathbb{R}^{3})$
is defined as the subset of distributions $f\in\mathcal{S}'(\mathbb{R}^{3})$ such that $\partial^{\beta}f\in \dot{B}^{s-k}_{p,r}(\mathbb{R}^{3})$
whenever $|\beta|=k$.
\end{itemize}
\end{definition}

\noindent\textbf{Remark 2.1}  Let $s\in \mathbb{R}$, $1\leq p,r\leq\infty$,  and $f\in\mathcal{S}'_{h}(\mathbb{R}^{3})$. Then $u\in
\dot{B}^{s}_{p,r}(\mathbb{R}^{3})$ if and only if there exists
$\{d_{j,r}\}_{j\in\mathbb{Z}}$ such that $d_{j,r}\ge0$,
$\|d_{j,r}\|_{\ell^{r}}=1$ and
$$
  \|\Delta_{j}u\|_{L^{p}}\lesssim
  d_{j,r}2^{-js}\|u\|_{\dot{B}^{s}_{p,r}} \ \ \text{for all }\
  j\in\mathbb{Z}.
$$

The fundamental idea in the proofs of Theorems \ref{th1.1} and \ref{th1.2} is to localize
system \eqref{eq1.3} through the Littlewood--Paley dyadic
decomposition, so we need the following definition of the Chemin--Lerner mixed time-space spaces, which was
first introduced by Chemin--Lerner \cite{CL95}.
\begin{definition}\label{de2.2}  For $0<T\leq\infty$, $s\in \mathbb{R}$ and
$1\leq p, r, \rho\leq\infty$. We define the mixed time-space $\mathcal{L}^{\rho}(0,T; \dot{B}^{s}_{p,r}(\mathbb{R}^{3}))$
as the completion of $\mathcal{C}([0,T]; \mathcal{S}(\mathbb{R}^{3}))$ by the norm
$$
  \|f\|_{\mathcal{L}^{\rho}_{T}(\dot{B}^{s}_{p,r})}:=\left(\sum_{j\in\mathbb{Z}}2^{srj}\left(\int_{0}^{T}
  \|\Delta_{j}f(\cdot,t)\|_{L^{p}}^{\rho}dt\right)^{\frac{r}{\rho}}\right)^{\frac{1}{r}}<\infty
$$
and with the standard modification for $\rho=\infty$ or $r=\infty$.
\end{definition}

\noindent\textbf{Remark 2.2} According to the  Minkowski inequality, it holds that
\begin{equation*}
  \|f\|_{\mathcal{L}^{\rho}_{T}(\dot{B}^{s}_{p,r})}\leq\|f\|_{L^{\rho}_{T}(\dot{B}^{s}_{p,r})} \ \ \  \text{if}\ \ \  \rho\leq r;\ \ \
  \|f\|_{L^{\rho}_{T}(\dot{B}^{s}_{p,r})}\leq \|f\|_{\mathcal{L}^{\rho}_{T}(\dot{B}^{s}_{p,r})} \ \ \ \text{if} \ \ \  r\leq \rho.
\end{equation*}
In particular, for $\rho=r=1$,  one has
\begin{equation*}
    \|f\|_{\mathcal{L}^{1}_{T}(\dot{B}^{s}_{p,1})}\approx\|f\|_{L^{1}_{T}(\dot{B}^{s}_{p,1})}.
\end{equation*}

In order to prove Theorem \ref{th1.2}, we need to introduce the following important weighted Chemin--Lerner type norm (see \cite{PZ11,PZ12}).

\begin{definition}\label{de2.3}  Let $f(t)\in L^{1}_{loc}(\mathbb{R}_{+})$, $f(t)\ge0$. For $1\leq p, r, \rho\leq\infty$, we define
$$
  \|u\|_{\mathcal{L}^{\rho}_{T,f}(\dot{B}^{s}_{p,r})}:=\left(\sum_{j\in\mathbb{Z}}2^{srj}\left(\int_{0}^{T}f(t)
  \|\Delta_{j}u(\cdot,t)\|_{L^{p}}^{\rho}dt\right)^{\frac{r}{\rho}}\right)^{\frac{1}{r}}<\infty
$$
and with the standard modification for $\rho=\infty$
or $r=\infty$.
\end{definition}

\subsection{Analytical tools in Besov spaces}

Let us recall the classical Bernstein inequality (see Lemma 2.1 in \cite{BCD11}).

\begin{lemma}\label{le2.4}
Let $\mathcal{B}$ be a ball, and $\mathcal{C}$ a ring in
$\mathbb{R}^{3}$. There exists a constant $C$ such that for any
positive real number $\lambda$, any nonnegative integer $k$ and any
couple of real numbers $(a,b)$ with $1\leq a\le b\leq \infty$, we
have
\begin{equation}\label{eq2.4}
   \operatorname{supp}\mathcal {F}(f)\subset\lambda\mathcal{B}\ \ \Rightarrow\ \   \sup_{|\alpha|=k}\|\partial^{\alpha}f\|_{L^{b}}\leq
   C^{k+1}\lambda^{k+3(\frac{1}{a}-\frac{1}{b})}\|f\|_{L^{a}},
\end{equation}
\begin{equation}\label{eq2.5}
   \operatorname{supp}\mathcal {F}(f)\subset\lambda\mathcal{C} \ \ \Rightarrow\ \   C^{-1-k}\lambda^{k}\|f\|_{L^{a}}\leq
   \sup_{|\alpha|=k}\|\partial^{\alpha}f\|_{L^{a}}\leq  C^{1+k}\lambda^{k}\|f\|_{L^{a}}.
\end{equation}
\end{lemma}

More generally, for any smooth homogeneous function $\sigma$ of degree $m$ on $\mathbb{R}^{3}/ \{0\}$ and $1\leq a\leq \infty$, it holds that
\begin{equation}\label{eq2.6}
  \operatorname{supp}\mathcal {F}(f)\subset\lambda\mathcal{C} \ \ \Rightarrow\ \  \|\sigma(D)f\|_{L^{a}}\lesssim
  \lambda^{k}\|f\|_{L^{a}}.
\end{equation}
An obvious consequence of \eqref{eq2.5} and \eqref{eq2.6} is that $\|\partial^{\alpha}f\|_{\dot{B}^{s}_{p,r}}\approx \|f\|_{\dot{B}^{s+k}_{p,r}}$ with multi-index $|\alpha|=k$ and $k\in\mathbb{N}$. Moreover,  the following lower bound for the integral involving the Laplacian
$-\Delta$, which can be regarded as a nonlinear generalization of \eqref{eq2.5},  will also be used,  for details,  see Lemma 8 in \cite{D10}.

\begin{lemma}\label{le2.5}
Suppose that $\operatorname{supp}\mathcal {F}(f)\subset \{\xi\in\mathbb{R}^{3}: \ \ K_{1}2^{j}\leq |\xi|\leq K_{2}2^{j}\}$ for some $K_1,K_2>0$ and  $j\in\mathbb{Z}$. Then there exists a constant $\kappa$  so that for all $1<p<\infty$, we have
\begin{equation}\label{eq2.7}
 -\int_{\mathbb{R}^{3}}\Delta f |f|^{p-2}f dx=(p-1)\int_{\mathbb{R}^{3}}|\nabla f|^{2}|f|^{p-2}dx\geq \kappa 2^{2j}\|f\|_{L^{p}}^{p},
\end{equation}
where $\kappa$ is a constant depending only on $p$, $K_1$ and $K_2$.
\end{lemma}

The following basic properties of Besov spaces are often used (see \cite{BCD11}).

\begin{lemma}\label{le5.2}  The following properties hold:
\begin{itemize}
\item [i)] Completeness: $\dot{B}^{s}_{p,r}(\mathbb{R}^{3})$ is a Banach space whenever $|s|<\frac{3}{p}$ or $s=\frac{3}{p}$ and $r=1$.

\item [ii)] Derivatives: There exists a universal constant $C$ such that
\begin{equation*}
    C^{-1}\|u\|_{\dot{B}^{s}_{p,r}}\leq \|\nabla u\|_{\dot{B}^{s-1}_{p,r}}\leq C\|u\|_{\dot{B}^{s}_{p,r}}.
\end{equation*}

\item [iii)] Fractional derivative: Let $\Lambda=\sqrt{-\Delta}$ and $\sigma\in\mathbb{R}$. Then the operator $\Lambda^{\sigma}$ is an isomorphism from $\dot{B}^{s}_{p,r}(\mathbb{R}^{3})$
to $\dot{B}^{s-\sigma}_{p,r}(\mathbb{R}^{3})$.

\item [iv)]  Imbedding: For $1\leq p_{1}\leq p_{2}\leq \infty$ and $1\leq r_{1}\leq r_{2}\leq \infty$, we have the continuous imbedding  $\dot{B}^{s}_{p_{1},r_{1}}(\mathbb{R}^{3})\hookrightarrow \dot{B}^{s-3(\frac{1}{p_{1}}-\frac{1}{p_{2}})}_{p_{2},r_{2}}(\mathbb{R}^{3})$.

\item [v)] Interpolation:  For  $s_{1},s_{2}\in \mathbb{R}$ such that $s_{1}<s_{2}$ and $\theta\in (0,1)$,  there exists a constant $C$
such that
\begin{align*}
    \|u\|_{\dot{B}^{s_{1}\theta+s_{2}(1-\theta)}_{p,r}}\leq C\|u\|_{\dot{B}^{s_{1}}_{p,r}}^{\theta}\|u\|_{\dot{B}^{s_{2}}_{p,r}}^{1-\theta}.
\end{align*}
\end{itemize}
\end{lemma}

The following crucial estimates for the product of two functions in the homogeneous Besov spaces are also used frequently throughout this paper (see Lemma 5.3 in \cite{ZL17}).

\begin{lemma}\label{le2.7}
Let $1\leq p_{1}, p_{2}\leq\infty$, $s_{1}\le \frac{3}{p_1}$,
$s_2\leq \min\{\frac{3}{p_{1}}, \frac{3}{p_2}\}$, and
$s_1+s_2>3\max\{0,\frac{1}{p_1}+\frac{1}{p_2}-1\}$. Assume that $f\in
\dot{B}^{s_1}_{p_1,1}(\mathbb{R}^{3})$,
$g\in\dot{B}^{s_2}_{p_2,1}(\mathbb{R}^{3})$. Then we have
$fg\in\dot{B}^{s_1+s_2-\frac{3}{p_1}}_{p_2,1}(\mathbb{R}^{3})$, and the following inequality holds:
\begin{equation}\label{eq2.8}
   \|fg\|_{\dot{B}^{s_1+s_2-\frac{3}{p_1}}_{p_2,1}}\lesssim\|f\|_{\dot{B}^{s_1}_{p_1,1}}\|g\|_{\dot{B}^{s_2}_{p_2,1}}.
\end{equation}
\end{lemma}

\subsection{Bilinear estimates}
 In this subsection, we recall the following bilinear estimates which are crucial steps to the proof of Theorem \ref{th1.1}, for the detailed proofs of these bilinear estimates, we refer the readers to  see \cite{ZL17,ZL22,ZZL15}. Here and in the sequel we denote $(d_{j})_{j\in\mathbb{Z}}$ a generic element of $l^{1}(\mathbb{Z})$ such that $d_{j}\ge0$ and
$\sum_{j\in\mathbb{Z}}d_{j}=1$.
\begin{lemma}\label{le2.8}
Let $1\leq p<\infty$. Then we have
\begin{align*}
   \|\Delta_{j}(u\cdot\nabla u)\|_{L^{1}_{t}(L^{p})}\lesssim 2^{j} \|\Delta_{j}(u\otimes u)\|_{L^{1}_{t}(L^{p})}\lesssim d_{j}2^{(1-\frac{3}{p})j} \|u\|_{\mathcal{L}^{\infty}_{T}(\dot{B}^{-1+\frac{3}{p}}_{p,1})}
   \|u\|_{L^{1}_{T}(\dot{B}^{1+\frac{3}{p}}_{p,1})},
\end{align*}
which implies that
\begin{equation}\label{eq2.9}
   \|u\cdot\nabla
   u\|_{L^{1}_{T}(\dot{B}^{-1+\frac{3}{p}}_{p,1})}=\|\nabla\cdot(u\otimes
   u)\|_{L^{1}_{T}(\dot{B}^{-1+\frac{3}{p}}_{p,1})}\lesssim\|u\|_{\mathcal{L}^{\infty}_{T}(\dot{B}^{-1+\frac{3}{p}}_{p,1})}
   \|u\|_{L^{1}_{T}(\dot{B}^{1+\frac{3}{p}}_{p,1})}.
\end{equation}
\end{lemma}

\begin{lemma}\label{le2.9}
Let $1\leq p, q<\infty$ and $\frac{1}{q}-\frac{1}{p}\geq-\min\{\frac{1}{3}, \frac{1}{2p}\}$.  Then we have
\begin{align*}
   \|\Delta_{j}(v\nabla(-\Delta)^{-1}v)\|_{L^{1}_{t}(L^{p})}\lesssim d_{j}2^{(1-\frac{3}{p})j} \|v\|_{L^{1}_{t}(\dot{B}^{\frac{3}{q}}_{q,1})}\|v\|_{\mathcal{L}^{\infty}_{t}(\dot{B}^{-2+\frac{3}{q}}_{q,1})},
\end{align*}
which implies that
\begin{align}\label{eq2.10}
   \|(v\nabla(-\Delta)^{-1}v)\|_{L^{1}_{t}(\dot{B}^{-1+\frac{3}{p}}_{p,1})}\lesssim \|v\|_{L^{1}_{t}(\dot{B}^{\frac{3}{q}}_{q,1})}\|v\|_{\mathcal{L}^{\infty}_{t}(\dot{B}^{-2+\frac{3}{q}}_{q,1})}.
\end{align}
\end{lemma}

\begin{lemma}\label{le2.10}
Let $1\leq p,q<\infty$ and $\frac{1}{p}+\frac{1}{q}>\frac{1}{3}$. Then we have
\begin{align*}
   \|\Delta_{j}(u\cdot\nabla v)\|_{L^{1}_{t}(L^{q})}\lesssim d_{j}2^{(2-\frac{3}{q})j} \big(\|u\|_{\mathcal{L}^{\infty}_{T}(\dot{B}^{-1+\frac{3}{p}}_{p,1})}\|v\|_{L^{1}_{t}(\dot{B}^{\frac{3}{q}}_{q,1})}+
   \|u\|_{L^{1}_{T}(\dot{B}^{1+\frac{3}{p}}_{p,1})}\|v\|_{\mathcal{L}^{\infty}_{t}(\dot{B}^{-2+\frac{3}{q}}_{q,1})}\big),
\end{align*}
which implies that
\begin{equation}\label{eq2.11}
   \|u\cdot\nabla v\|_{L^{1}_{T}(\dot{B}^{-2+\frac{3}{q}}_{q,1})}\lesssim\|u\|_{\mathcal{L}^{\infty}_{T}(\dot{B}^{-1+\frac{3}{p}}_{p,1})}\|v\|_{L^{1}_{t}(\dot{B}^{\frac{3}{q}}_{q,1})}+
   \|u\|_{L^{1}_{T}(\dot{B}^{1+\frac{3}{p}}_{p,1})}\|v\|_{\mathcal{L}^{\infty}_{t}(\dot{B}^{-2+\frac{3}{q}}_{q,1})}.
\end{equation}
\end{lemma}

\begin{lemma}\label{le2.11}
Let $1\leq q,r<\infty$ and $\frac{1}{q}-\frac{1}{r}<\frac{1}{3}<\frac{1}{q}+\frac{1}{r}$. Then we have
\begin{align*}
   \|\Delta_{j}(w\nabla(-\Delta)^{-1}v)\|_{L^{1}_{t}(L^{q})}\lesssim d_{j}2^{(1-\frac{3}{q})j} \big(\|w\|_{\mathcal{L}^{\infty}_{T}(\dot{B}^{-2+\frac{3}{r}}_{r,1})}\|v\|_{L^{1}_{t}(\dot{B}^{\frac{3}{q}}_{q,1})}+
   \|w\|_{L^{1}_{T}(\dot{B}^{\frac{3}{r}}_{r,1})}\|v\|_{\mathcal{L}^{\infty}_{t}(\dot{B}^{-2+\frac{3}{q}}_{q,1})}\big),
\end{align*}
which implies that
\begin{equation}\label{eq2.12}
   \|w\nabla(-\Delta)^{-1}v\|_{L^{1}_{T}(\dot{B}^{-1+\frac{3}{q}}_{q,1})}\lesssim\|w\|_{\mathcal{L}^{\infty}_{T}(\dot{B}^{-2+\frac{3}{r}}_{r,1})}\|v\|_{L^{1}_{t}(\dot{B}^{\frac{3}{q}}_{q,1})}+
   \|w\|_{L^{1}_{T}(\dot{B}^{\frac{3}{r}}_{r,1})}\|v\|_{\mathcal{L}^{\infty}_{t}(\dot{B}^{-2+\frac{3}{q}}_{q,1})}.
\end{equation}
\end{lemma}

\begin{lemma}\label{le2.12}
Let $1\leq p,r<\infty$ and $\frac{1}{p}+\frac{1}{r}>\frac{1}{3}$. Then we have
\begin{align*}
   \|\Delta_{j}(u\cdot\nabla w)\|_{L^{1}_{t}(L^{r})}\lesssim d_{j}2^{(2-\frac{3}{r})j} \big(\|u\|_{\mathcal{L}^{\infty}_{T}(\dot{B}^{-1+\frac{3}{p}}_{p,1})}\|w\|_{L^{1}_{t}(\dot{B}^{\frac{3}{r}}_{r,1})}+
   \|u\|_{L^{1}_{T}(\dot{B}^{1+\frac{3}{p}}_{p,1})}\|w\|_{\mathcal{L}^{\infty}_{t}(\dot{B}^{-2+\frac{3}{r}}_{r,1})}\big),
\end{align*}
which implies that
\begin{equation}\label{eq2.13}
   \|u\cdot\nabla v\|_{L^{1}_{T}(\dot{B}^{-2+\frac{3}{q}}_{q,1})}\lesssim\|u\|_{\mathcal{L}^{\infty}_{T}(\dot{B}^{-1+\frac{3}{p}}_{p,1})}\|w\|_{L^{1}_{t}(\dot{B}^{\frac{3}{r}}_{r,1})}+
   \|u\|_{L^{1}_{T}(\dot{B}^{1+\frac{3}{p}}_{p,1})}\|w\|_{\mathcal{L}^{\infty}_{t}(\dot{B}^{-2+\frac{3}{r}}_{r,1})}.
\end{equation}
\end{lemma}

\begin{lemma}\label{le2.13}
Let $1\leq q, r<\infty$ and $\frac{1}{r}-\frac{1}{q}<\frac{1}{3}$.  Then we have
\begin{align*}
   \|\Delta_{j}(v\nabla(-\Delta)^{-1}v)\|_{L^{1}_{t}(L^{r})}\lesssim d_{j}2^{(1-\frac{3}{r})j} \|v\|_{L^{1}_{t}(\dot{B}^{\frac{3}{q}}_{q,1})}\|v\|_{\mathcal{L}^{\infty}_{t}(\dot{B}^{-2+\frac{3}{q}}_{q,1})},
\end{align*}
which implies that
\begin{align}\label{eq2.14}
   \|(v\nabla(-\Delta)^{-1}v)\|_{L^{1}_{t}(\dot{B}^{-1+\frac{3}{r}}_{r,1})}\lesssim \|v\|_{L^{1}_{t}(\dot{B}^{\frac{3}{q}}_{q,1})}\|v\|_{\mathcal{L}^{\infty}_{t}(\dot{B}^{-2+\frac{3}{q}}_{q,1})}.
\end{align}
\end{lemma}

\section{Estimates of the velocity field $u$}

The purpose of this section is to derive the desired estimates for the horizontal components $u^{h}=(u^{1},u^{2})$ and the vertical component $u^{3}$ of the velocity field in the framework of weighted Chemin--Lerner type spaces. The main idea is that we introduce some weighted functions and weighted norms to eliminate the difficulties caused by the nonlinear terms involving $u^{3}$ and $w$. Thus we set
\begin{equation}\label{eq3.1}
  f_{1}(t):=\|u^{3}(\cdot,t)\|_{\dot{B}^{1+\frac{3}{p}}_{p,1}},  \ \
  f_{2}(t):=\|u^{3}(\cdot,t)\|_{\dot{B}^{\frac{3}{p}}_{p,1}}^{2},  \ \ f_{3}(t):=\|w(\cdot,t)\|_{\dot{B}^{\frac{3}{r}}_{r,1}}.
\end{equation}
For three positive real numbers $\lambda_{1}$, $\lambda_{2}$ and $\lambda_{3}$, we denote $\vec{\lambda}=(\lambda_{1}, \lambda_{2}, \lambda_{3})$, and introduce the following three weighted functions:
\begin{align*}
  &u^{h}_{\vec{\lambda}}:=u\exp\big\{-\lambda_{1}\int_{0}^{t}f_{1}(\tau)d\tau-\lambda_{2}\int_{0}^{t}f_{2}(\tau)d\tau-\lambda_{3}\int_{0}^{t}f_{3}(\tau)d\tau\big\},\\
  &\pi_{\vec{\lambda}}:=\pi\exp\big\{-\lambda_{1}\int_{0}^{t}f_{1}(\tau)d\tau-\lambda_{2}\int_{0}^{t}f_{2}(\tau)d\tau-\lambda_{3}\int_{0}^{t}f_{3}(\tau)d\tau\big\},\\
 & v_{\vec{\lambda}}:=v\exp\big\{-\lambda_{1}\int_{0}^{t}f_{1}(\tau)d\tau-\lambda_{2}\int_{0}^{t}f_{2}(\tau)d\tau-\lambda_{3}\int_{0}^{t}f_{3}(\tau)d\tau\big\}.
\end{align*}

\subsection{Estimates of the pressure $\pi$}

Notice that, using the divergence free condition $\nabla\cdot u=0$, the term $\nabla\cdot(u\cdot\nabla u)$ has a nice structure:
\begin{equation*}
  \nabla\cdot(u\cdot\nabla u)=\operatorname{div}_{h}\operatorname{div}_{h}(u^{h}\otimes
  u^{h})+2\partial_{3}\operatorname{div}_{h}(u^{3}u^{h})+\partial_{3}^{2}(u^{3})^{2},
\end{equation*}
where for a vector field $u=(u^{1},u^{2},u^{3})=(u^{h},u^{3})$, we
denote $\operatorname{div}_{h}u^{h}=\partial_{1}u^{1}+\partial_{2}u^{2}$.  Thus, by taking the divergence $\operatorname{div}$ to the first equations of \eqref{eq1.3}
yields that
\begin{equation}\label{eq3.2}
  -\Delta\pi=\operatorname{div}_{h}\operatorname{div}_{h}(u^{h}\otimes
  u^{h})+2\partial_{3}\operatorname{div}_{h}(u^{3}u^{h})+\partial_{3}^{2}(u^{3})^{2}+\nabla\cdot(v\nabla(-\Delta)^{-1}v).
\end{equation}
Multiplying \eqref{eq3.2} by the weighted function $\exp\big\{-\sum_{i=1}^{3}\lambda_{i}\int_{0}^{t}f_{i}(\tau)d\tau\big\}$ and applying the divergence free condition $\partial_{3}u^{3}=-\operatorname{div}_{h}u^{h}$, we arrive at
\begin{align}\label{eq3.3}
  \nabla\pi_{\vec{\lambda}}=\nabla(-\Delta)^{-1}\Big[\operatorname{div}_{h}\operatorname{div}_{h}(u^{h}\otimes u^{h}_{\vec{\lambda}})
  +2\partial_{3}\operatorname{div}_{h}(u^{3}u^{h}_{\vec{\lambda}})
  -2\partial_{3}(u^{3}\operatorname{div}_{h}u^{h}_{\vec{\lambda}})
  +\nabla\cdot(v\nabla(-\Delta)^{-1}v_{\vec{\lambda}})\Big].
\end{align}
Applying the dyadic operator $\Delta_{j}$ to \eqref{eq3.3}, then taking $L^{1}_{t}(L^{p})$-norm and using the
H\"{o}lder inequality and Bernstein inequality \eqref{eq2.4} yield that
\begin{align}\label{eq3.4}
  \|\Delta_{j}(\nabla\pi_{\vec{\lambda}})\|_{L^{1}_{t}(L^{p})}&
  \lesssim2^{j}\big(\|\Delta_{j}(u^{h}\otimes u^{h}_{\vec{\lambda}})\|_{L^{1}_{t}(L^{p})}
  +\|\Delta_{j}(u^{3}u^{h}_{\vec{\lambda}})\|_{L^{1}_{t}(L^{p})}\big)\nonumber\\
  &+\|\Delta_{j}(u^{3}\operatorname{div}_{h}u^{h}_{\vec{\lambda}})\|_{L^{1}_{t}(L^{p})}
  +\|\Delta_{j}(v\nabla(-\Delta)^{-1}v_{\vec{\lambda}})\|_{L^{1}_{t}(L^{p})}.
\end{align}
The first three terms on the right-hand side of \eqref{eq3.4} have been estimated in \cite{PZ11, PZ12},  and the last term has been bounded in Lemma \ref{le2.9},  thus we obtain
\begin{align}\label{eq3.5}
  \|\Delta_{j}(u^{h}\otimes
  u^{h}_{\vec{\lambda}})\|_{L^{1}_{t}(L^{p})}\lesssim
  d_{j}2^{-\frac{3j}{p}}\|u^{h}\|_{\mathcal{L}^{\infty}_{t}(\dot{B}^{-1+\frac{3}{p}}_{p,1})}
  \|u^{h}_{\vec{\lambda}}\|_{L^{1}_{t}(\dot{B}^{1+\frac{3}{p}}_{p,1})},
\end{align}
\begin{equation}\label{eq3.6}
  \|\Delta_{j}(u^{3}u^{h}_{\vec{\lambda}})\|_{L^{1}_{t}(L^{p})}\lesssim
  d_{j}2^{-\frac{3j}{p}}\big(\|u^{h}_{\vec{\lambda}}\|_{L^{1}_{t}(\dot{B}^{1+\frac{3}{p}}_{p,1})}^{\frac{1}{2}}
  \|u^{h}_{\vec{\lambda}}\|_{\mathcal{L}^{1}_{t,f_{2}}(\dot{B}^{-1+\frac{3}{p}}_{p,1})}^{\frac{1}{2}}
  +\|u^{h}_{\vec{\lambda}}\|_{\mathcal{L}^{1}_{t,f_{1}}(\dot{B}^{-1+\frac{3}{p}}_{p,1})}\big),
\end{equation}
\begin{equation}\label{eq3.7}
  \|\Delta_{j}(u^{3}\operatorname{div}_{h}u^{h}_{\vec{\lambda}})\|_{L^{1}_{t}(L^{p})}\lesssim
  d_{j}2^{j(1-\frac{3}{p})}\big(\|u^{h}_{\vec{\lambda}}\|_{L^{1}_{t}(\dot{B}^{1+\frac{3}{p}}_{p,1})}^{\frac{1}{2}}
  \|u^{h}_{\vec{\lambda}}\|_{\mathcal{L}^{1}_{t,f_{2}}(\dot{B}^{-1+\frac{3}{p}}_{p,1})}^{\frac{1}{2}}
  +\|u^{h}_{\vec{\lambda}}\|_{\mathcal{L}^{1}_{t,f_{1}}(\dot{B}^{-1+\frac{3}{p}}_{p,1})}\big),
\end{equation}
\begin{equation}\label{eq3.8}
   \|\Delta_{j}(v\nabla(-\Delta)^{-1}v_{\vec{\lambda}})\|_{L^{1}_{t}(L^{p})}\lesssim d_{j}2^{(1-\frac{3}{p})j} \|v\|_{\mathcal{L}^{\infty}_{t}(\dot{B}^{-2+\frac{3}{q}}_{q,1})}\|v_{\vec{\lambda}}\|_{L^{1}_{t}(\dot{B}^{\frac{3}{q}}_{q,1})}.
\end{equation}
Taking the above estimates \eqref{eq3.5}-\eqref{eq3.8} into \eqref{eq3.4}, we get
\begin{align}\label{eq3.9}
  \|\Delta_{j}(\nabla\pi_{\vec{\lambda}})&\|_{L^{1}_{t}(L^{p})}\lesssim
  d_{j}2^{j(1-\frac{3}{p})}\Big[\|u^{h}_{\vec{\lambda}}\|_{L^{1}_{t}(\dot{B}^{1+\frac{3}{p}}_{p,1})}^{\frac{1}{2}}
  \|u^{h}_{\vec{\lambda}}\|_{\mathcal{L}^{1}_{t,f_{2}}(\dot{B}^{-1+\frac{3}{p}}_{p,1})}^{\frac{1}{2}}
  +\|u^{h}_{\vec{\lambda}}\|_{\mathcal{L}^{1}_{t,f_{1}}(\dot{B}^{-1+\frac{3}{p}}_{p,1})}\nonumber\\
  &+\|u^{h}\|_{\mathcal{L}^{\infty}_{t}(\dot{B}^{-1+\frac{3}{p}}_{p,1})}
  \|u^{h}_{\vec{\lambda}}\|_{L^{1}_{t}(\dot{B}^{1+\frac{3}{p}}_{p,1})}
  +\|v\|_{\mathcal{L}^{\infty}_{t}(\dot{B}^{-2+\frac{3}{q}}_{q,1})}
  \|v_{\vec{\lambda}}\|_{L^{1}_{t}(\dot{B}^{\frac{3}{q}}_{q,1})}\Big],
\end{align}
which clearly  implies that
\begin{align}\label{eq3.10}
  \|\nabla\pi_{\vec{\lambda}}&\|_{L^{1}_{t}(\dot{B}^{-1+\frac{3}{p}}_{p,1})}
  \lesssim \|u^{h}_{\vec{\lambda}}\|_{L^{1}_{t}(\dot{B}^{1+\frac{3}{p}}_{p,1})}^{\frac{1}{2}}
  \|u^{h}_{\vec{\lambda}}\|_{\mathcal{L}^{1}_{t,f_{2}}(\dot{B}^{-1+\frac{3}{p}}_{p,1})}^{\frac{1}{2}}
  +\|u^{h}_{\vec{\lambda}}\|_{\mathcal{L}^{1}_{t,f_{1}}(\dot{B}^{-1+\frac{3}{p}}_{p,1})}\nonumber
  \\
  &+\|u^{h}\|_{\mathcal{L}^{\infty}_{t}(\dot{B}^{-1+\frac{3}{p}}_{p,1})}
  \|u^{h}_{\vec{\lambda}}\|_{L^{1}_{t}(\dot{B}^{1+\frac{3}{p}}_{p,1})}+\|v\|_{\mathcal{L}^{\infty}_{t}(\dot{B}^{-2+\frac{3}{q}}_{q,1})}
  \|v_{\vec{\lambda}}\|_{L^{1}_{t}(\dot{B}^{\frac{3}{q}}_{q,1})}.
\end{align}

On the other hand, in order to deal with the vertical component $u^{3}$, we also need the
following two estimates from \cite{PZ12}:
\begin{equation}\label{eq3.11}
  \|\Delta_{j}(u^{3}u^{h})\|_{L^{1}_{t}(L^{p})}\lesssim
  d_{j}2^{-\frac{3j}{p}}\big(\|u^{h}\|_{\mathcal{L}^{\infty}_{t}(\dot{B}^{-1+\frac{3}{p}}_{p,1})}
  \|u^{3}\|_{L^{1}_{t}(\dot{B}^{1+\frac{3}{p}}_{p,1})}
  +\|u^{h}\|_{L^{1}_{t}(\dot{B}^{1+\frac{3}{p}}_{p,1})}
  \|u^{3}\|_{\mathcal{L}^{\infty}_{t}(\dot{B}^{-1+\frac{3}{p}}_{p,1})}\big),
\end{equation}
\begin{equation}\label{eq3.12}
  \|\Delta_{j}(u^{3}\operatorname{div}_{h}u^{h})\|_{L^{1}_{t}(L^{p})}\lesssim
  d_{j}2^{j(1-\frac{3}{p})}\big(\|u^{h}\|_{\mathcal{L}^{\infty}_{t}(\dot{B}^{-1+\frac{3}{p}}_{p,1})}
  \|u^{3}\|_{L^{1}_{t}(\dot{B}^{1+\frac{3}{p}}_{p,1})}
  \!+\!\|u^{h}\|_{L^{1}_{t}(\dot{B}^{1+\frac{3}{p}}_{p,1})}
  \|u^{3}\|_{\mathcal{L}^{\infty}_{t}(\dot{B}^{-1+\frac{3}{p}}_{p,1})}\big).
\end{equation}
Based on the above two estimates, we can exactly follow the same lines as derivation of \eqref{eq3.10} by taking $\lambda_{1}=\lambda_{2}=\lambda_{3}=0$ to obtain that the pressure $\pi$ satisfies the following estimate:
\begin{align}\label{eq3.13}
  \|\nabla\pi\|_{L^{1}_{t}(\dot{B}^{-1+\frac{3}{p}}_{p,1})}&\le
  C\big(\|u^{h}\|_{\mathcal{L}^{\infty}_{t}(\dot{B}^{-1+\frac{3}{p}}_{p,1})}
  \big(\|u^{h}\|_{L^{1}_{t}(\dot{B}^{1+\frac{3}{p}}_{p,1})}+\|u^{3}\|_{L^{1}_{t}(\dot{B}^{1+\frac{3}{p}}_{p,1})}\big)\nonumber\\
  &+\|u^{h}\|_{L^{1}_{t}(\dot{B}^{1+\frac{3}{p}}_{p,1})}\|u^{3}\|_{\mathcal{L}^{\infty}_{t}(\dot{B}^{-1+\frac{3}{p}}_{p,1})}
  +\|v\|_{\mathcal{L}^{\infty}_{t}(\dot{B}^{-2+\frac{3}{q}}_{q,1})}
  \|v\|_{L^{1}_{t}(\dot{B}^{\frac{3}{q}}_{q,1})}\big).
\end{align}

\subsection{Estimate of the horizontal components $u^{h}$}
Considering the first equations of \eqref{eq1.3}, it is clear that the horizontal components $u^{h}$ satisfies the following equations:
\begin{equation}\label{eq3.14}
  \partial_{t}u^{h}_{\vec{\lambda}}+(\sum_{i=1}^{3}\lambda_{i}f_{i}(t))u^{h}_{\vec{\lambda}}-\Delta
  u^{h}_{\vec{\lambda}}=-u\cdot\nabla
  u^{h}_{\vec{\lambda}}-\nabla_{h}\pi_{\vec{\lambda}}-v\nabla_{h}(-\Delta)^{-1}v_{\vec{\lambda}}.
\end{equation}
Applying the operator $\Delta_{j}$ to \eqref{eq3.14}, then taking $L^{2}$
inner product of the resulting equations with
$|\Delta_{j}u^{h}_{\vec{\lambda}}|^{p-2}\Delta_{j}u^{h}_{\vec{\lambda}}$, we obtain
\begin{align}\label{eq3.15}
  \frac{1}{p}\frac{d}{dt}\|\Delta_{j}u^{h}_{\vec{\lambda}}\|_{L^{p}}^{p}&
  +(\sum_{i=1}^{3}\lambda_{i}f_{i}(t))\|\Delta_{j}u^{h}_{\vec{\lambda}}\|_{L^{p}}^{p}
 -\int_{\mathbb{R}^{3}}\Delta\Delta_{j}u^{h}_{\vec{\lambda}}\cdot|\Delta_{j}u^{h}_{\vec{\lambda}}|^{p-2}\Delta_{j}u^{h}_{\vec{\lambda}}dx\nonumber\\
  &=
  -\int_{\mathbb{R}^{3}}\Delta_{j}\big(u\cdot\nabla u^{h}_{\vec{\lambda}}
  +\nabla_{h}\pi_{\vec{\lambda}}+v\nabla_{h}(-\Delta)^{-1}v_{\vec{\lambda}}\big)
  |\Delta_{j}u^{h}_{\vec{\lambda}}|^{p-2}\Delta_{j}u^{h}_{\vec{\lambda}}dx.
\end{align}
Thanks to Lemma \ref{le2.5}, there exists a positive constant $\kappa$ such  that
\begin{equation*}
  -\int_{\mathbb{R}^{3}}\Delta\Delta_{j}u^{h}_{\vec{\lambda}}\cdot|\Delta_{j}u^{h}_{\vec{\lambda}}|^{p-2}\Delta_{j}u^{h}_{\vec{\lambda}}dx\geq
  \kappa2^{2j}\|\Delta_{j}u^{h}_{\vec{\lambda}}\|_{L^{p}}^{p},
\end{equation*}
whence a similar argument as that in \cite{D01} gives rise to
\begin{align*}
  \frac{d}{dt}\|\Delta_{j}u^{h}_{\vec{\lambda}}\|_{L^{p}}&+(\sum_{i=1}^{3}\lambda_{i}f_{i}(t))\|\Delta_{j}u^{h}_{\vec{\lambda}}\|_{L^{p}}
  +\kappa2^{2j}\|\Delta_{j}u^{h}_{\vec{\lambda}}\|_{L^{p}}\nonumber\\
  &\leq \|\Delta_{j}(u\cdot\nabla u^{h}_{\vec{\lambda}})\|_{L^{p}}+\|\Delta_{j}
  \nabla_{h}\pi_{\vec{\lambda}}\|_{L^{p}}+\|\Delta_{j}(v\nabla_{h}(-\Delta)^{-1}v_{\vec{\lambda}})\|_{L^{p}}.
\end{align*}
Integrating the above inequality on $[0,t]$ yields that
\begin{align}\label{eq3.16}
  \|\Delta_{j}u^{h}_{\vec{\lambda}}\|_{L^{\infty}_{t}(L^{p})}&+(\sum_{i=1}^{3}\lambda_{i}f_{i}(t))\|\Delta_{j}u^{h}_{\vec{\lambda}}\|_{L^{1}_{t}(L^{p})}
  +\kappa2^{2j}\|\Delta_{j}u^{h}_{\vec{\lambda}}\|_{L^{1}_{t}(L^{p})}\leq 2^{j}\big(\|\Delta_{j}(u^{h}\otimes  u^{h}_{\vec{\lambda}})\|_{L^{1}_{t}(L^{p})}\nonumber\\
  &+\|\Delta_{j}(u^{3}u^{h}_{\vec{\lambda}})\|_{L^{1}_{t}(L^{p})}\big)
  +\|\Delta_{j}
  \nabla_{h}\pi_{\vec{\lambda}}\|_{L^{1}_{t}(L^{p})}+\|\Delta_{j}(v\nabla_{h}(-\Delta)^{-1}v_{\vec{\lambda}})\|_{L^{1}_{t}(L^{p})}.
\end{align}
The right-hand side of \eqref{eq3.16} has been estimated in \eqref{eq3.5}, \eqref{eq3.6},  \eqref{eq3.8} and \eqref{eq3.9}, thus there exists a positive constant $C_{1}$ such that
\begin{align}\label{eq3.17}
  \|u^{h}_{\vec{\lambda}}\|_{\mathcal{L}^{\infty}_{t}(\dot{B}^{-1+\frac{3}{p}}_{p,1})}
  &+\sum_{i=1}^{3}\lambda_{i}\|u^{h}_{\vec{\lambda}}\|_{\mathcal{L}^{1}_{t,f_{i}}(\dot{B}^{-1+\frac{3}{p}}_{p,1})}
  +\kappa\|u^{h}_{\vec{\lambda}}\|_{L^{1}_{t}(\dot{B}^{1+\frac{3}{p}}_{p,1})}\nonumber\\
  &\leq \|u_{0}^{h}\|_{\dot{B}^{-1+\frac{3}{p}}_{p,1}}+\frac{\kappa}{4}\|u^{h}_{\vec{\lambda}}\|_{L^{1}_{t}(\dot{B}^{1+\frac{3}{p}}_{p,1})}
  +C_{1}\big(\|u^{h}\|_{\mathcal{L}^{\infty}_{t}(\dot{B}^{-1+\frac{3}{p}}_{p,1})}
  \|u^{h}_{\vec{\lambda}}\|_{L^{1}_{t}(\dot{B}^{1+\frac{3}{p}}_{p,1})}\nonumber\\
  &+\|u^{h}_{\vec{\lambda}}\|_{\mathcal{L}^{1}_{t,f_{1}}(\dot{B}^{-1+\frac{3}{p}}_{p,1})}
  +\|u^{h}_{\vec{\lambda}}\|_{\mathcal{L}^{1}_{t,f_{2}}(\dot{B}^{-1+\frac{3}{p}}_{p,1})}
  +\|v\|_{\mathcal{L}^{\infty}_{t}(\dot{B}^{-2+\frac{3}{q}}_{q,1})}
  \|v_{\vec{\lambda}}\|_{L^{1}_{t}(\dot{B}^{\frac{3}{q}}_{q,1})}
  \big).
\end{align}

\subsection{Estimate of the vertical component $u^{3}$}
Observing that the vertical component $u^{3}$ satisfies the following equation:
\begin{equation}\label{eq3.18}
  \partial_{t}u^{3}-\Delta
  u^{3}=-u\cdot\nabla
  u^{3}-\partial_{3}\pi-v\partial_{3}(-\Delta)^{-1}v.
\end{equation}
As the derivation of \eqref{eq3.16}, and using the following identity:
$$
  u\cdot\nabla u^{3}=\operatorname{div}(uu^{3})=\operatorname{div}_{h}(u^{h}u^{3})-2u^{3}\operatorname{div}_{h}u^{h},
$$
one has
\begin{align}\label{eq3.19}
  \|\Delta_{j}u^{3}\|_{L^{\infty}_{t}(L^{p})}+\kappa2^{2j}\|\Delta_{j}u^{3}\|_{L^{1}_{t}(L^{p})}&\leq \|u_{0}^{3}\|_{L^{p}}
  +C\big(\|\Delta_{j}(u\cdot\nabla u^{3})\|_{L^{1}_{t}(L^{p})}\nonumber\\
  &+\|\Delta_{j}\partial_{3}\pi\|_{L^{1}_{t}(L^{p})}
  +\|\Delta_{j}(v\partial_{3}(-\Delta)^{-1}v)\|_{L^{1}_{t}(L^{p})}\big)\nonumber\\
  &\leq \|u_{0}^{3}\|_{L^{p}}
  +C\big(2^{j}\|\Delta_{j}(u^{h}u^{3})\|_{L^{1}_{t}(L^{p})}+\|\Delta_{j}(u^{3}\operatorname{div}_{h}u^{h})\|_{L^{1}_{t}(L^{p})}\nonumber\\
  &+\|\Delta_{j}\partial_{3}\pi\|_{L^{1}_{t}(L^{p})}
  +\|\Delta_{j}(v\partial_{3}(-\Delta)^{-1}v)\|_{L^{1}_{t}(L^{p})}\big).
\end{align}
The right-hand side of \eqref{eq3.19} has been estimated in \eqref{eq3.11}, \eqref{eq3.12},  \eqref{eq3.13} and \eqref{eq2.10}, and substituting these estimates into \eqref{eq3.19}, we obtain that there exists a positive constant $C_{2}$ such that
\begin{align}\label{eq3.20}
  \|u^{3}\|_{\mathcal{L}^{\infty}_{t}(\dot{B}^{-1+\frac{3}{p}}_{p,1})}
  &+\kappa\|u^{3}\|_{L^{1}_{t}(\dot{B}^{1+\frac{3}{p}}_{p,1})}
  \leq \|u_{0}^{3}\|_{\dot{B}^{-1+\frac{3}{p}}_{p,1}}
  +C_{2}\big(\|u^{h}\|_{\mathcal{L}^{\infty}_{t}(\dot{B}^{-1+\frac{3}{p}}_{p,1})}
  (\|u^{h}\|_{L^{1}_{t}(\dot{B}^{1+\frac{3}{p}}_{p,1})}+\|u^{3}\|_{L^{1}_{t}(\dot{B}^{1+\frac{3}{p}}_{p,1})})\nonumber\\
  &+\|u^{h}\|_{L^{1}_{t}(\dot{B}^{1+\frac{3}{p}}_{p,1})}\|u^{3}\|_{\mathcal{L}^{\infty}_{t}(\dot{B}^{-1+\frac{3}{p}}_{p,1})}
  +\|v\|_{\mathcal{L}^{\infty}_{t}(\dot{B}^{-2+\frac{3}{q}}_{q,1})}
  \|v\|_{L^{1}_{t}(\dot{B}^{\frac{3}{q}}_{q,1})}\big).
\end{align}

\section{Estimates of the densities $v$ and $w$}

In this section, we intend to derive the estimates for the charged densities $v$ and $w$. As we pointed out before, the crucial ingredient is to introduce the proper weighted functions to eliminate the difficulties caused by the nonlinear terms of the third and fourth equations of system \eqref{eq1.3}, and we shall use different weighted functions to tackle with $v$ and $w$.

\subsection{Estimate of  density v}
To deal with $v$, we mainly use $f_{1}(t)$ to eliminate the difficulties caused the term $u\cdot\nabla v$, and the weighted function $f_{3}(t)$ to eliminate the difficulties caused the term $\nabla\cdot(w\nabla(-\Delta)^{-1}v)$.  It follows the third equation of \eqref{eq1.3} that
\begin{equation}\label{eq4.1}
  \partial_{t}v_{\vec{\lambda}}+(\sum_{i=1}^{3}\lambda_{i}f_{i}(t))v_{\vec{\lambda}}-\Delta v_{\vec{\lambda}}=-u\cdot\nabla
  v_{\vec{\lambda}}+\nabla\cdot(w\nabla(-\Delta)^{-1}v_{\vec{\lambda}}).
\end{equation}
Applying the dyadic operator $\Delta_{j}$ to \eqref{eq4.1}, then taking
$L^{2}$ inner product of the resulting equation with
$|\Delta_{j}v_{\vec{\lambda}}|^{q-2}\Delta_{j}v_{\vec{\lambda}}$ and applying Lemma \ref{le2.5},  we see that
\begin{align}\label{eq4.2}
  \frac{1}{q}\frac{d}{dt}\|\Delta_{j}v_{\vec{\lambda}}\|_{L^{q}}^{q}&+(\sum_{i=1}^{3}\lambda_{i}f_{i}(t))\|\Delta_{j}v_{\vec{\lambda}}\|_{L^{q}}^{q}+
  \kappa2^{2j}\|\Delta_{j}v_{\vec{\lambda}}\|_{L^{q}}^{q}\nonumber\\
  &\leq
  -\int_{\mathbb{R}^{3}}\big(\Delta_{j}(u\cdot\nabla v_{\vec{\lambda}})+\Delta_{j}
  \nabla\cdot(w\nabla(-\Delta)^{-1}v_{\vec{\lambda}}\big)|\Delta_{j}v_{\vec{\lambda}}|^{q-2}\Delta_{j}v_{\vec{\lambda}}dx.
\end{align}
Moreover, applying Bony's paraproduct decomposition \eqref{eq2.3}, one has
\begin{equation*}
  u\cdot\nabla v_{\vec{\lambda}}=T_{u}\nabla v_{\vec{\lambda}}+T_{\nabla
  v_{\vec{\lambda}}}u+R(u,\nabla v_{\vec{\lambda}}),
\end{equation*}
which combining the standard commutator's argument gives us to
\begin{align}\label{eq4.3}
 \int_{\mathbb{R}^{3}}\Delta_{j}(T_{u}\nabla
 v_{\vec{\lambda}})|\Delta_{j}v_{\vec{\lambda}}|^{q-2}\Delta_{j}v_{\vec{\lambda}}dx&=\sum_{|j'-j|\leq
 5}\int_{\mathbb{R}^{3}}[\Delta_{j};S_{j'-1}u]\Delta_{j'}\nabla
 v_{\vec{\lambda}}|\Delta_{j}v_{\vec{\lambda}}|^{q-2}\Delta_{j}v_{\vec{\lambda}}dx\nonumber\\
 &+\sum_{|j'-j|\leq5}\int_{\mathbb{R}^{3}}(S_{j'-1}u-S_{j-1}u)\Delta_{j}\Delta_{j'}\nabla
 v_{\vec{\lambda}}|\Delta_{j}v_{\vec{\lambda}}|^{q-2}\Delta_{j}v_{\vec{\lambda}}dx\nonumber\\
 &-\frac{1}{q}\int_{\mathbb{R}^{3}}S_{j-1}(\operatorname{div}
 u)\Delta_{j}\Delta_{j'}
 v_{\vec{\lambda}}|\Delta_{j}v_{\vec{\lambda}}|^{q-2}\Delta_{j}v_{\vec{\lambda}}dx.
\end{align}
Hence, taking the above estimate \eqref{eq4.3} into \eqref{eq4.2},  and using the divergence free condition $\nabla\cdot u=0$ and the
argument for the $L^{q}$ energy estimate in \cite{D01}, we obtain
\begin{align}\label{eq4.4}
  \|\Delta_{j}v_{\vec{\lambda}}\|_{L^{q}}&+\int_{0}^{t}(\sum_{i=1}^{3}\lambda_{i}f_{i}(\tau))\|\Delta_{j}v_{\vec{\lambda}}(\tau)\|_{L^{q}}d\tau
  +\kappa2^{2j}\int_{0}^{t}\|\Delta_{j}v_{\vec{\lambda}}(\tau)\|_{L^{q}}d\tau\leq
  \|\Delta_{j}v_{0}\|_{L^{q}}\nonumber\\
  &+C\big(\sum_{|j'-j|\leq 5}\big(\|[\Delta_{j};S_{j'-1}u]\Delta_{j'}\nabla
  v_{\vec{\lambda}}\|_{L^{1}_{t}(L^{q})}+\|(S_{j'-1}u-S_{j-1}u)\Delta_{j}\Delta_{j'}\nabla
  v_{\vec{\lambda}}\|_{L^{1}_{t}(L^{q})}\big)\nonumber\\
  &+\|\Delta_{j}(T_{\nabla
  v_{\vec{\lambda}}}u)\|_{L^{1}_{t}(L^{q})}+\|\Delta_{j}(R(u,\nabla
  v_{\vec{\lambda}}))\|_{L^{1}_{t}(L^{q})}+\|\Delta_{j}
  \nabla\cdot(w\nabla(-\Delta)^{-1}v_{\vec{\lambda}})\|_{L^{1}_{t}(L^{q})}\big).
\end{align}
In the following we estimate the terms on the right-hand side of \eqref{eq4.4} one by one. Applying Lemma \ref{le2.4} and Definition \ref{de2.1}, the first two terms can be estimated as
\begin{align}\label{eq4.5}
  \sum_{|j'-j|\leq 5}&\|[\Delta_{j};S_{j'-1}u]\Delta_{j'}\nabla
  v_{\vec{\lambda}}\|_{L^{1}_{t}(L^{q})}\nonumber\\
  &\lesssim \sum_{|j'-j|\leq
  5}\big(\|S_{j'-1}\nabla u^{h}\|_{L^{1}_{t}(L^{\infty})}\|\Delta_{j'}v_{\vec{\lambda}}\|_{L^{\infty}_{t}(L^{q})}+\int_{0}^{t}\|S_{j'-1}\nabla
  u^{3}(\tau)\|_{L^{\infty}}\|\Delta_{j'}v_{\vec{\lambda}}(\tau)\|_{L^{q}}d\tau\big)\nonumber\\
  &\lesssim \sum_{|j'-j|\leq 5}\big(d_{j'}2^{j'(2-\frac{3}{q})}\|u^{h}\|_{L^{1}_{t}(\dot{B}^{1+\frac{3}{p}}_{p,1})}
  \|v_{\vec{\lambda}}\|_{\mathcal{L}^{\infty}_{t}(\dot{B}^{-2+\frac{3}{q}}_{q,1})}
  +\int_{0}^{t}\|u^{3}(\tau)\|_{\dot{B}^{1+\frac{3}{p}}_{p,1}}
  \|\Delta_{j'}v_{\vec{\lambda}}(\tau)\|_{L^{q}}d\tau\big)\nonumber\\
  &\lesssim d_{j}2^{j(2-\frac{3}{q})}\big(\|u^{h}\|_{L^{1}_{t}(\dot{B}^{1+\frac{3}{p}}_{p,1})}
  \|v_{\vec{\lambda}}\|_{\mathcal{L}^{\infty}_{t}(\dot{B}^{-2+\frac{3}{q}}_{q,1})}
  +\|v_{\vec{\lambda}}\|_{\mathcal{L}^{1}_{t,f_{1}}(\dot{B}^{-2+\frac{3}{q}}_{q,1})}\big),
\end{align}
\begin{align}\label{eq4.6}
  \sum_{|j'-j|\leq 5}&\|(S_{j'-1}u-S_{j-1}u)\Delta_{j}\Delta_{j'}\nabla
  v_{\vec{\lambda}}\|_{L^{1}_{t}(L^{q})}\nonumber\\
  &\lesssim \sum_{|j'-j|\leq
  5}\big(\|(S_{j'-1}\nabla u^{h}-S_{j-1}\nabla u^{h})\|_{L^{1}_{t}(L^{\infty})}\|\Delta_{j}v_{\vec{\lambda}}\|_{L^{\infty}_{t}(L^{q})}\nonumber\\
  &\ \ \ \  +\int_{0}^{t}\|(S_{j'-1}\nabla u^{3}-S_{j-1}\nabla u^{3})(\tau)\|_{L^{\infty}}\|\Delta_{j}v_{\vec{\lambda}}(\tau)\|_{L^{q}}d\tau\big)\nonumber\\
  &\lesssim d_{j}2^{j(2-\frac{3}{q})}\|u^{h}\|_{L^{1}_{t}(\dot{B}^{1+\frac{3}{p}}_{p,1})}
  \|v_{\vec{\lambda}}\|_{\mathcal{L}^{\infty}_{t}(\dot{B}^{-2+\frac{3}{q}}_{q,1})}
  +\int_{0}^{t}\|u^{3}(\tau)\|_{\dot{B}^{1+\frac{3}{p}}_{p,1}}
  \|\Delta_{j}v_{\vec{\lambda}}(\tau)\|_{L^{q}}d\tau\nonumber\\
  &\lesssim d_{j}2^{j(2-\frac{3}{q})}\big(\|u^{h}\|_{L^{1}_{t}(\dot{B}^{1+\frac{3}{p}}_{p,1})}
  \|v_{\vec{\lambda}}\|_{\mathcal{L}^{\infty}_{t}(\dot{B}^{-2+\frac{3}{q}}_{q,1})}
  +\|v_{\vec{\lambda}}\|_{\mathcal{L}^{1}_{t,f_{1}}(\dot{B}^{-2+\frac{3}{q}}_{q,1})}\big).
\end{align}
For the term involving $T_{\nabla
  v_{\vec{\lambda}}}u$, we consider the following two cases: in the case $1<q\leq p<6$, one can choose $\widetilde{q}$ ($1< \widetilde{q}\leq
\infty$) such that $\frac{1}{q}=\frac{1}{p}+\frac{1}{\widetilde{q}}$, then applying
Lemma \ref{le2.4} yields that
\begin{align}\label{eq4.7}
  \|\Delta_{j}(T_{\nabla v_{\vec{\lambda}}}u)\|_{L^{1}_{t}(L^{q})}&\lesssim
  \sum_{|j'-j|\leq5}\big(\|S_{j'-1}\nabla_{h}v_{\vec{\lambda}}\|_{L^{\infty}_{t}(L^{\widetilde{q}})}\|\Delta_{j'}u^{h}\|_{L^{1}_{t}(L^{p})}
  \nonumber\\
  &\ \ \ \  +\int_{0}^{t}\|S_{j'-1}\partial_{3}v_{\vec{\lambda}}(\tau)
  \|_{L^{\widetilde{q}}}\|\Delta_{j'}u^{3}(\tau)\|_{L^{p}}d\tau\big)\nonumber\\
  &\lesssim \sum_{k\leq
   j'-2}2^{k(1+\frac{3}{q}-\frac{3}{\widetilde{q}})}\|\Delta_{k}v_{\vec{\lambda}}\|_{L^{\infty}_{t}(L^{q})}\|\Delta_{j'}u^{h}\|_{L^{1}_{t}(L^{p})}\nonumber\\
   &\ \ \ \  +2^{-j(1+\frac{3}{p})}\sum_{|j'-j|\leq5}\sum_{k\leq
   j'-2}2^{k(1+\frac{3}{p})}\int_{0}^{t}\|u^{3}(\tau)\|_{\dot{B}^{1+\frac{3}{p}}_{p,1}}\|\Delta_{k}v_{\vec{\lambda}}(\tau)\|_{L^{q}}d\tau\nonumber\\
   &\lesssim \sum_{k\leq
   j'-2}d_{k}2^{k(3+\frac{3}{p}-\frac{3}{q})}\|v_{\vec{\lambda}}\|_{\mathcal{L}^{\infty}_{t}(\dot{B}^{-2+\frac{3}{q}}_{q,1})}\|\Delta_{j'}u^{h}\|_{L^{1}_{t}(L^{p})}\nonumber\\
   &\ \ \ \  +2^{-j(1+\frac{3}{p})}\sum_{|j'-j|\leq5}\sum_{k\leq
   j'-2}2^{k(3+\frac{3}{p}-\frac{3}{q})}d_{k}\|v_{\vec{\lambda}}\|_{\mathcal{L}^{1}_{t,f_{1}}(\dot{B}^{-2+\frac{3}{q}}_{q,1})}\nonumber\\
   &\lesssim d_{j}2^{j(2-\frac{3}{q})}\big(\|u^{h}\|_{L^{1}_{t}(\dot{B}^{1+\frac{3}{p}}_{p,1})}
  \|v_{\vec{\lambda}}\|_{\mathcal{L}^{\infty}_{t}(\dot{B}^{-2+\frac{3}{q}}_{q,1})}
  +\|v_{\vec{\lambda}}\|_{\mathcal{L}^{1}_{t,f_{1}}(\dot{B}^{-2+\frac{3}{q}}_{q,1})}\big),
\end{align}
while in the case $1< p<q$,  one can directly estimate that
\begin{align}\label{eq4.8}
  \|\Delta_{j}(T_{\nabla v_{\vec{\lambda}}}u)\|_{L^{1}_{t}(L^{q})}&\lesssim 2^{(\frac{3}{p}-\frac{3}{q})j}
  \sum_{|j'-j|\leq5}\big(\|S_{j'-1}\nabla_{h}v_{\vec{\lambda}}\|_{L^{\infty}_{t}(L^{\infty})}\|\Delta_{j'}u^{h}\|_{L^{1}_{t}(L^{p})}
  \nonumber\\
  &\ \ \ \  +\int_{0}^{t}\|S_{j'-1}\partial_{3}v_{\vec{\lambda}}(\tau)
  \|_{L^{\infty}}\|\Delta_{j'}u^{3}(\tau)\|_{L^{p}}d\tau\big)\nonumber\\
  &\lesssim 2^{(\frac{3}{p}-\frac{3}{q})j}\sum_{k\leq
   j'-2}2^{(1+\frac{3}{q})k}\|\Delta_{k}v_{\vec{\lambda}}\|_{L^{\infty}_{t}(L^{q})}\|\Delta_{j'}u^{h}\|_{L^{1}_{t}(L^{p})}\nonumber\\
   &\ \ \ \  +d_{j}2^{-j(1+\frac{3}{q})}\sum_{|j'-j|\leq5}\sum_{k\leq
   j'-2}2^{(1+\frac{3}{q})k}\int_{0}^{t}\|u^{3}(\tau)\|_{\dot{B}^{1+\frac{3}{p}}_{p,1}}\|\Delta_{k}v_{\vec{\lambda}}(\tau)\|_{L^{q}}d\tau\nonumber\\
   &\lesssim 2^{(\frac{3}{p}-\frac{3}{q})j}\sum_{k\leq
   j'-2}d_{k}2^{3k}\|v_{\vec{\lambda}}\|_{\mathcal{L}^{\infty}_{t}(\dot{B}^{-2+\frac{3}{q}}_{q,1})}\|\Delta_{j'}u^{h}\|_{L^{1}_{t}(L^{p})}\nonumber\\
   &\ \ \ \  + 2^{-j(1+\frac{3}{q})}\sum_{|j'-j|\leq5}\sum_{k\leq
   j'-2}2^{3k}d_{k}\|v_{\vec{\lambda}}\|_{\mathcal{L}^{1}_{t,f_{1}}(\dot{B}^{-2+\frac{3}{q}}_{q,1})}\nonumber\\
   &\lesssim d_{j}2^{j(2-\frac{3}{q})}\big(\|u^{h}\|_{L^{1}_{t}(\dot{B}^{1+\frac{3}{p}}_{p,1})}
  \|v_{\vec{\lambda}}\|_{\mathcal{L}^{\infty}_{t}(\dot{B}^{-2+\frac{3}{q}}_{q,1})}
  +\|v_{\vec{\lambda}}\|_{\mathcal{L}^{1}_{t,f_{1}}(\dot{B}^{-2+\frac{3}{q}}_{q,1})}\big).
\end{align}
For the remaining term $\|\Delta_{j}(R(u,\nabla
  v_{\vec{\lambda}}))\|_{L^{1}_{t}(L^{q})}$, we split the estimate into the following two parts: In the case $\frac{1}{3}<\frac{1}{p}+\frac{1}{q}\leq1$, we get
\begin{align}\label{eq4.9}
  \|\Delta_{j}(R(u,\nabla
  v_{\vec{\lambda}}))\|_{L^{1}_{t}(L^{q})}&\lesssim
  2^{(1+\frac{3}{p})j}\sum_{j'\geq
  j-N_{0}}\|\Delta_{j'}u^{h}\|_{L^{1}_{t}(L^{p})}\|\widetilde{\Delta}_{j'}v_{\vec{\lambda}}\|_{L^{\infty}_{t}(L^{q})}\nonumber\\
  &\ \ \ \ + 2^{(1+\frac{3}{p})j}\sum_{j'\geq
  j-N_{0}}\int_{0}^{t}\|\widetilde{\Delta}_{j'}u^{3}(\tau)\|_{L^{p}}\|\Delta_{j'}v_{\vec{\lambda}}(\tau)\|_{L^{q}}d\tau\nonumber\\
  &\lesssim 2^{(1+\frac{3}{p})j}\sum_{j'\geq
  j-N_{0}}d_{j'}2^{j'(1-\frac{3}{p}-\frac{3}{q})}\|u^{h}\|_{L^{1}_{t}(\dot{B}^{1+\frac{3}{p}}_{p,1})}
  \|v_{\vec{\lambda}}\|_{\mathcal{L}^{\infty}_{t}(\dot{B}^{-2+\frac{3}{q}}_{q,1})}\nonumber\\
  &\ \ \ \ +2^{(1+\frac{3}{p})j}\sum_{j'\geq
  j-N_{0}}2^{-(1+\frac{3}{p})j'}\int_{0}^{t}\|u^{3}(\tau)\|_{\dot{B}^{1+\frac{3}{p}}_{p,1}}\|\Delta_{j'}v_{\vec{\lambda}}(\tau)\|_{L^{q}}d\tau\nonumber\\
  &\lesssim d_{j}2^{j(2-\frac{3}{q})}\|u^{h}\|_{L^{1}_{t}(\dot{B}^{1+\frac{3}{p}}_{p,1})}\big(
  \|v_{\vec{\lambda}}\|_{\mathcal{L}^{\infty}_{t}(\dot{B}^{-2+\frac{3}{q}}_{q,1})}
  +\|v_{\vec{\lambda}}\|_{\mathcal{L}^{1}_{t,f_{1}}(\dot{B}^{-2+\frac{3}{q}}_{q,1})}\big);
\end{align}
while in the case $\frac{1}{p}+\frac{1}{q}>1$, we choose  $q'$ ($1<q'<\infty$)  such that $\frac{1}{q}+\frac{1}{q'}=1$, then applying
Lemma \ref{le2.4} again yields that
\begin{align}\label{eq4.10}
 \|\Delta_{j}(R(u,\nabla
  v_{\vec{\lambda}}))\|_{L^{1}_{t}(L^{q})}&\lesssim
  2^{j(3-\frac{3}{q})}\sum_{j'\geq
  j-N_{0}}\|\Delta_{j'}u^{h}\|_{L^{1}_{t}(L^{q'})}\|\widetilde{\Delta}_{j'}\nabla_{h}v_{\vec{\lambda}}\|_{L^{\infty}_{t}(L^{q})}\nonumber\\
  &\ \ \ \ +2^{j(3-\frac{3}{q})}\sum_{j'\geq
  j-N_{0}}\int_{0}^{t}\|\widetilde{\Delta}_{j'}u^{3}(\tau)\|_{L^{q'}}\|\Delta_{j'}\partial_{3}v_{\vec{\lambda}}(\tau)\|_{L^{q}}d\tau\nonumber\\
  &\lesssim 2^{j(3-\frac{3}{q})}\sum_{j'\geq
  j-N_{0}}2^{3j'(\frac{1}{p}-\frac{1}{q'})}\|\Delta_{j'}u^{h}\|_{L^{1}_{t}(L^{p})}\|\widetilde{\Delta}_{j'}\nabla_{h}v_{\vec{\lambda}}\|_{L^{\infty}_{t}(L^{q})}\nonumber\\
  &\ \ \ \ +2^{j(3-\frac{3}{q})}\sum_{j'\geq
  j-N_{0}}2^{3j'(\frac{1}{p}-\frac{1}{q'})}\int_{0}^{t}\|\widetilde{\Delta}_{j'}u^{3}(\tau)\|_{L^{p}}
  \|\Delta_{j'}\partial_{3}v_{\vec{\lambda}}(\tau)\|_{L^{q}}d\tau\nonumber\\
  &\lesssim 2^{j(3-\frac{3}{q})}\sum_{j'\geq
  j-N_{0}}d_{j'}2^{-j'}\big(\|u^{h}\|_{L^{1}_{t}(\dot{B}^{1+\frac{3}{p}}_{p,1})}
  \|v_{\vec{\lambda}}\|_{\mathcal{L}^{\infty}_{t}(\dot{B}^{-2+\frac{3}{q}}_{q,1})}
  +\|v_{\vec{\lambda}}\|_{\mathcal{L}^{1}_{t,f_{1}}(\dot{B}^{-2+\frac{3}{q}}_{q,1})}\big)\nonumber\\
  &\lesssim d_{j}2^{j(2-\frac{3}{q})}\big(\|u^{h}\|_{L^{1}_{t}(\dot{B}^{1+\frac{3}{p}}_{p,1})}
  \|v_{\vec{\lambda}}\|_{\mathcal{L}^{\infty}_{t}(\dot{B}^{-2+\frac{3}{q}}_{q,1})}+
  \|v_{\vec{\lambda}}\|_{\mathcal{L}^{1}_{t,f_{1}}(\dot{B}^{-2+\frac{3}{q}}_{q,1})}\big).
\end{align}
Putting all above estimates together, we conclude that
\begin{align*}
  &\|\Delta_{j}v_{\vec{\lambda}}\|_{L^{q}}+\sum_{i=1}^{3}\lambda_{i}\int_{0}^{t}f_{i}(\tau)\|\Delta_{j}v_{\vec{\lambda}}(\tau)\|_{L^{q}}d\tau
    +\kappa2^{2j}\int_{0}^{t}\|\Delta_{j}v_{\vec{\lambda}}(\tau)\|_{L^{q}}d\tau\leq
  \|\Delta_{j}v_{0}\|_{L^{q}}\nonumber\\
  &+Cd_{j}2^{j(2-\frac{3}{q})}\big(\|u^{h}\|_{L^{1}_{t}(\dot{B}^{1+\frac{3}{p}}_{p,1})}
  \|v_{\vec{\lambda}}\|_{\mathcal{L}^{\infty}_{t}(\dot{B}^{-2+\frac{3}{q}}_{q,1})}
  +\|v_{\vec{\lambda}}\|_{\mathcal{L}^{1}_{t,f_{1}}(\dot{B}^{-2+\frac{3}{q}}_{q,1})}\big)+C \|\Delta_{j} \nabla\cdot(w\nabla(-\Delta)^{-1}v_{\vec{\lambda}})\|_{L^{1}_{t}(L^{q})},
\end{align*}
which directly implies that
\begin{align}\label{eq4.11}
  &\|v_{\vec{\lambda}}\|_{\mathcal{L}^{\infty}_{t}(\dot{B}^{-2+\frac{3}{q}}_{q,1})}
  +\sum_{i=1}^{3}\lambda_{i}\|v_{\vec{\lambda}}\|_{\mathcal{L}^{1}_{t,f_{i}}(\dot{B}^{-2+\frac{3}{q}}_{q,1})}
  +\kappa\|v_{\vec{\lambda}}\|_{L^{1}_{t}(\dot{B}^{\frac{3}{q}}_{q,1})}\leq
  \|v_{0}\|_{\dot{B}^{-2+\frac{3}{q}}_{q,1}}\nonumber\\
  &+C\big(\|u^{h}\|_{L^{1}_{t}(\dot{B}^{1+\frac{3}{p}}_{p,1})}
  \|v_{\vec{\lambda}}\|_{\mathcal{L}^{\infty}_{t}(\dot{B}^{-2+\frac{3}{q}}_{q,1})}
  +\|v_{\vec{\lambda}}\|_{\mathcal{L}^{1}_{t,f_{1}}(\dot{B}^{-2+\frac{3}{q}}_{q,1})}
  +\|w\nabla(-\Delta)^{-1}v_{\vec{\lambda}}\|_{\mathcal{L}^{1}_{t}(\dot{B}^{-1+\frac{3}{q}}_{q,1})}\big).
\end{align}
According to the Minkowski's inequality, it is readily to see that
\begin{align*}
  \|w\nabla(-\Delta)^{-1}v_{\vec{\lambda}}\|_{\mathcal{L}^{1}_{t}(\dot{B}^{-1+\frac{3}{q}}_{q,1})}\approx
  \int_{0}^{t}\|w(\tau)\nabla(-\Delta)^{-1}v_{\vec{\lambda}}(\tau)\|_{\dot{B}^{-1+\frac{3}{q}}_{q,1}}d\tau.
\end{align*}
Then we can apply Lemma \ref{le2.7} by setting $s_{1}=\frac{3}{r}$, $s_{2}=-1+\frac{3}{q}$,  $p_{1}=r$ and $p_{2}=q$
 to obtain that
\begin{align*}
    \|w\nabla(-\Delta)^{-1}v_{\vec{\lambda}}\|_{\dot{B}^{-1+\frac{3}{q}}_{q,1}}
    \lesssim \|w\|_{\dot{B}^{\frac{3}{r}}_{r,1}}\|\nabla(-\Delta)^{-1}v_{\vec{\lambda}}\|_{\dot{B}^{-1+\frac{3}{q}}_{q,1}}
     \lesssim \|w\|_{\dot{B}^{\frac{3}{r}}_{r,1}}\|v_{\vec{\lambda}}\|_{\dot{B}^{-2+\frac{3}{q}}_{q,1}},
\end{align*}
which yields that
\begin{align}\label{eq4.12}
  \|w\nabla(-\Delta)^{-1}v_{\vec{\lambda}}\|_{\mathcal{L}^{1}_{t}(\dot{B}^{-1+\frac{3}{q}}_{q,1})}\lesssim \int_{0}^{t}\|w(\tau)\|_{\dot{B}^{\frac{3}{r}}_{r,1}}\|v_{\vec{\lambda}}(\tau)\|_{\dot{B}^{-2+\frac{3}{q}}_{q,1}}d\tau
  \lesssim \|v_{\vec{\lambda}}\|_{\mathcal{L}^{1}_{t,f_{3}}(\dot{B}^{-2+\frac{3}{q}}_{q,1})}.
\end{align}
Taking \eqref{eq4.12} into \eqref{eq4.11}, we obtain that there exists a positive constant $C_{3}$ such that
\begin{align}\label{eq4.13}
  \|v_{\vec{\lambda}}\|_{\mathcal{L}^{\infty}_{t}(\dot{B}^{-2+\frac{3}{q}}_{q,1})}&
  +\sum_{i=1}^{3}\lambda_{i}\|v_{\vec{\lambda}}\|_{\mathcal{L}^{1}_{t,f_{i}}(\dot{B}^{-2+\frac{3}{q}}_{q,1})}
  +\kappa\|v_{\vec{\lambda}}\|_{L^{1}_{t}(\dot{B}^{\frac{3}{q}}_{q,1})}\leq
  \|v_{0}\|_{\dot{B}^{-2+\frac{3}{q}}_{q,1}}\nonumber\\
  &+C_{3}\big(\|u^{h}\|_{L^{1}_{t}(\dot{B}^{1+\frac{3}{p}}_{p,1})}
  \|v_{\vec{\lambda}}\|_{\mathcal{L}^{\infty}_{t}(\dot{B}^{-2+\frac{3}{q}}_{q,1})}
  +\|v_{\vec{\lambda}}\|_{\mathcal{L}^{1}_{t,f_{1}}(\dot{B}^{-2+\frac{3}{q}}_{q,1})}
  +\|v_{\vec{\lambda}}\|_{\mathcal{L}^{1}_{t,f_{3}}(\dot{B}^{-2+\frac{3}{q}}_{q,1})}\big).
\end{align}

\subsection{Estimate of density $w$}
For any positive real number $\lambda_{1}$, recall that $f_{1}(t)=\|u^{3}(t)\|_{\dot{B}^{1+\frac{3}{p}}_{p,1}}$,  and we denote
\begin{equation*}
w_{\lambda_{1}}:=w\exp(-\lambda_{1}\int_{0}^{t}f_{1}(\tau)d\tau),\ \ \
v_{\lambda_{1}}:=v\exp(-\lambda_{1}\int_{0}^{t}f_{1}(\tau)d\tau).
\end{equation*}
It follows the fourth equation of \eqref{eq1.1} that
\begin{equation}\label{eq4.14}
  \partial_{t}w_{\lambda_{1}}+\lambda_{1}f_{1}(t)w_{\lambda_{1}}-\Delta w_{\lambda_{1}}=-u\cdot\nabla
  w_{\lambda_{1}}+\nabla\cdot(v\nabla(-\Delta)^{-1}v_{\lambda_{1}}).
\end{equation}
 Arguing like the derivation of \eqref{eq4.1} yields that
 \begin{align}\label{eq4.15}
  &\|w_{\lambda_{1}}\|_{\mathcal{L}^{\infty}_{t}(\dot{B}^{-2+\frac{3}{r}}_{r,1})}
  +\lambda_{1}\|w_{\lambda_{1}}\|_{\mathcal{L}^{1}_{t,f_{1}}(\dot{B}^{-2+\frac{3}{r}}_{r,1})}
  +\kappa\|w_{\lambda_{1}}\|_{L^{1}_{t}(\dot{B}^{\frac{3}{r}}_{r,1})}\leq
  \|w_{0}\|_{\dot{B}^{-2+\frac{3}{r}}_{r,1}}\nonumber\\
  &+C\big(\|u^{h}\|_{L^{1}_{t}(\dot{B}^{1+\frac{3}{p}}_{p,1})}
  \|w_{\lambda_{1}}\|_{\mathcal{L}^{\infty}_{t}(\dot{B}^{-2+\frac{3}{r}}_{r,1})}
  +\|w_{\lambda_{1}}\|_{\mathcal{L}^{1}_{t,f_{1}}(\dot{B}^{-2+\frac{3}{r}}_{r,1})}
  +\|v\nabla(-\Delta)^{-1}v_{\lambda_{1}}\|_{L^{1}_{t}(\dot{B}^{-1+\frac{3}{r}}_{r,1})}\big).
\end{align}
Applying Lemma \ref{le2.13}, one has
\begin{equation*}
\|v\nabla(-\Delta)^{-1}v_{\lambda_{1}}\|_{L^{1}_{t}(\dot{B}^{-1+\frac{3}{r}}_{r,1})}\lesssim \|v_{\lambda_{1}}\|_{L^{1}_{t}(\dot{B}^{\frac{3}{q}}_{q,1})}\|v\|_{\mathcal{L}^{\infty}_{t}(\dot{B}^{-2+\frac{3}{q}}_{q,1})},
\end{equation*}
which back to \eqref{eq4.15}, we conclude that  there exists a positive constant $C_{4}$ such that
\begin{align}\label{eq4.16}
  &\|w_{\lambda_{1}}\|_{\mathcal{L}^{\infty}_{t}(\dot{B}^{-2+\frac{3}{r}}_{r,1})}
  +\lambda_{1}\|w_{\lambda_{1}}\|_{\mathcal{L}^{1}_{t,f_{1}}(\dot{B}^{-2+\frac{3}{r}}_{r,1})}
  +\kappa\|w_{\lambda_{1}}\|_{L^{1}_{t}(\dot{B}^{\frac{3}{r}}_{r,1})}\leq
  \|w_{0}\|_{\dot{B}^{-2+\frac{3}{r}}_{r,1}}\nonumber\\
  &+C_{4}\big(\|u^{h}\|_{L^{1}_{t}(\dot{B}^{1+\frac{3}{p}}_{p,1})}
  \|w_{\lambda_{1}}\|_{\mathcal{L}^{\infty}_{t}(\dot{B}^{-2+\frac{3}{r}}_{r,1})}
  +\|w_{\lambda_{1}}\|_{\mathcal{L}^{1}_{t,f_{1}}(\dot{B}^{-2+\frac{3}{r}}_{r,1})}+\|v\|_{\mathcal{L}^{\infty}_{t}(\dot{B}^{-2+\frac{3}{q}}_{q,1})}
  \|v_{\lambda_{1}}\|_{L^{1}_{t}(\dot{B}^{\frac{3}{q}}_{q,1})}\big).
\end{align}

\section{Proofs of Theorems \ref{th1.1} and \ref{th1.2}}

The proof of Theorem \ref{th1.1} is simple. Once one gets the above desired bilinear estimates \eqref{eq2.9}--\eqref{eq2.14},  one can follow exactly the same procedure as \cite{ZZL15} to prove that there exists $T>0$ such that the system \eqref{eq1.3}  admits a unique
solution $(u, v,w)$ on $[0,T]$ satisfying \eqref{eq1.8}. Moreover, if the initial data is sufficiently small, then the above local solution is actually a global one,
for details, please see \cite{ZZL15}.

Now we present the proof of Theorem \ref{th1.2}. Let us denote by $T_{*}$ the maximal existence time of  local solution $(u, v,w)$ satisfying \eqref{eq1.8}. Then to prove Theorem \ref{th1.2},
it suffices to prove $T_{*}=\infty$ under the initial condition \eqref{eq1.9}. We prove it by contradiction. Assume that $T_{*}<\infty$, based on the estimates \eqref{eq3.17}, \eqref{eq3.20},
\eqref{eq4.13} and \eqref{eq4.16}, let $\eta$ be a small enough positive constant which the exact value will be determined later, we define $T_{\eta}$ by
\begin{align}\label{eq5.1}
  T_{\eta}:=\max\Big\{t\in[0,T_{*}):\ \|u^{h}\|_{\mathcal{L}^{\infty}_{t}(\dot{B}^{-1+\frac{3}{p}}_{p,1})}
  +\kappa\|u^{h}\|_{L^{1}_{t}(\dot{B}^{1+\frac{3}{p}}_{p,1})}+\|v\|_{\mathcal{L}^{\infty}_{t}(\dot{B}^{-2+\frac{3}{q}}_{q,1})}
  +\kappa\|v\|_{L^{1}_{t}(\dot{B}^{\frac{3}{q}}_{q,1})}\leq
  \eta \Big\}.
\end{align}
Taking $\lambda_{1}\geq 2C_{1}$, $\lambda_{2}\geq 2C_{1}$ and $\eta\leq
\frac{\kappa}{4C_{1}}$, we can derive from \eqref{eq3.17} to get that
\begin{align}\label{eq5.2}
   \|u^{h}_{\vec{\lambda}}\|_{\mathcal{L}^{\infty}_{t}(\dot{B}^{-1+\frac{3}{p}}_{p,1})}
   +\frac{\kappa}{2}\|u^{h}_{\vec{\lambda}}\|_{L^{1}_{t}(\dot{B}^{1+\frac{3}{p}}_{p,1})}
   &\leq \|u_{0}^{h}\|_{\dot{B}^{-1+\frac{3}{p}}_{p,1}}+\frac{\kappa}{4}\|v_{\vec{\lambda}}\|_{L^{1}_{t}(\dot{B}^{\frac{3}{q}}_{q,1})}.
\end{align}
On the other hand, taking $\lambda_{1}\geq 2C_{3}$,  $\lambda_{3}\geq 2C_{3}$ in  \eqref{eq4.13},  and $\eta\leq
 \frac{\kappa}{2C_{3}}$ one obtains that
\begin{align}\label{eq5.3}
   \|v_{\vec{\lambda}}\|_{\mathcal{L}^{\infty}_{t}(\dot{B}^{-2+\frac{3}{q}}_{q,1})}
    +2\kappa\|v_{\vec{\lambda}}\|_{L^{1}_{t}(\dot{B}^{\frac{3}{q}}_{q,1})}\leq
  2\|v_{0}\|_{\dot{B}^{-2+\frac{3}{q}}_{q,1}}.
\end{align}
As a consequence,  we obtain from \eqref{eq5.2}--\eqref{eq5.3} that for all $t\leq T_{\eta}$, it holds that
\begin{align}\label{eq5.4}
   &\|u^{h}\|_{\mathcal{L}^{\infty}_{t}(\dot{B}^{-1+\frac{3}{p}}_{p,1})}
   +\kappa\|u^{h}\|_{L^{1}_{t}(\dot{B}^{1+\frac{3}{p}}_{p,1})}
   +\|v\|_{\mathcal{L}^{\infty}_{t}(\dot{B}^{-2+\frac{3}{q}}_{q,1})}
   +\kappa\|v\|_{L^{1}_{t}(\dot{B}^{\frac{3}{q}}_{q,1})}\nonumber\\
   &\leq2\big(\|u_{0}^{h}\|_{\dot{B}^{-1+\frac{3}{p}}_{p,1}}+ \|v_{0}\|_{\dot{B}^{-2+\frac{3}{q}}_{q,1}}\big)
   \times \exp\Big\{\int_{0}^{t}(\lambda_{1}f_{1}(\tau)+\lambda_{2}f_{2}(\tau)+\lambda_{3}f_{3}(\tau))d\tau\Big\}
   \nonumber\\
   &=2\big(\|u_{0}^{h}\|_{\dot{B}^{-1+\frac{3}{p}}_{p,1}}+ \|v_{0}\|_{\dot{B}^{-2+\frac{3}{q}}_{q,1}}\big)
   \times \exp\Big\{\int_{0}^{t}(\lambda_{1}\|u^{3}(\tau)\|_{\dot{B}^{1+\frac{3}{p}}_{p,1}}+\lambda_{2}\|u^{3}(\tau)\|_{\dot{B}^{\frac{3}{p}}_{p,1}}^{2}
   +\lambda_{3}\|w(\tau)\|_{\dot{B}^{\frac{3}{r}}_{r,1}})d\tau\Big\}.
\end{align}
Thanks to \eqref{eq3.20}, by choosing $\eta\leq
\ \frac{\kappa}{2C_{2}}$, it holds that for all $t\leq T_{\eta}$,
\begin{align}\label{eq5.5}
  \|u^{3}\|_{\mathcal{L}^{\infty}_{t}(\dot{B}^{-1+\frac{3}{p}}_{p,1})}
  +\kappa\|u^{3}\|_{L^{1}_{t}(\dot{B}^{1+\frac{3}{p}}_{p,1})}
  &\leq
  2\|u_{0}^{3}\|_{\dot{B}^{-1+\frac{3}{p}}_{p,1}}+2\eta.
\end{align}
Back to \eqref{eq4.16}, by taking $\lambda_{1}>2C_{4}$ and $\eta\leq
\ \frac{\kappa}{2C_{4}}$, one gets
\begin{align*}
  \|w_{\lambda_{1}}\|_{\mathcal{L}^{\infty}_{t}(\dot{B}^{-2+\frac{3}{r}}_{r,1})}
    +\kappa\|w_{\lambda_{1}}\|_{L^{1}_{t}(\dot{B}^{\frac{3}{r}}_{r,1})}\leq
  2\|w_{0}\|_{\dot{B}^{-2+\frac{3}{r}}_{r,1}}+\eta,
\end{align*}
which using \eqref{eq5.5} yields that
\begin{align}\label{eq5.6}
  \|w\|_{\mathcal{L}^{\infty}_{t}(\dot{B}^{-2+\frac{3}{r}}_{r,1})}
    +\kappa\|w\|_{L^{1}_{t}(\dot{B}^{\frac{3}{r}}_{r,1})}&\leq
  \big(2\|w_{0}\|_{\dot{B}^{-2+\frac{3}{r}}_{r,1}}+\eta\big)\times \exp\Big\{\int_{0}^{t}\lambda_{1}\|u^{3}(\tau)\|_{\dot{B}^{1+\frac{3}{p}}_{p,1}}d\tau\Big\}\nonumber\\
  &\leq \big(2\|w_{0}\|_{\dot{B}^{-2+\frac{3}{r}}_{r,1}}+\eta\big) \times \exp\Big\{\frac{2\lambda_{1}}{\kappa}(\|u_{0}^{3}\|_{\dot{B}^{-1+\frac{3}{p}}_{p,1}}+\eta) \Big\}.
\end{align}
Besides, it follows the interpolation inequality in Lemma \ref{eq2.6} that
\begin{align}\label{eq5.7}
  \|u^{3}\|_{L^{2}_{t}(\dot{B}^{\frac{3}{p}}_{p,1})}^{2}&\leq C
  \|u^{3}\|_{L^{\infty}_{t}(\dot{B}^{-1+\frac{3}{p}}_{p,1})}
  \|u^{3}\|_{L^{1}_{t}(\dot{B}^{1+\frac{3}{p}}_{p,1})}\nonumber\\
  &\leq
  C\|u^{3}\|_{\mathcal{L}^{\infty}_{t}(\dot{B}^{-1+\frac{3}{p}}_{p,1})}
  \|u^{3}\|_{L^{1}_{t}(\dot{B}^{1+\frac{3}{p}}_{p,1})}\nonumber\\
  &\leq \frac{C}{\kappa}(\|u_{0}^{3}\|_{\dot{B}^{-1+\frac{3}{p}}_{p,1}}+\eta)^{2}.
\end{align}
Taking above estimates \eqref{eq5.5}--\eqref{eq5.7} into \eqref{eq5.4}, we obtain that there exists a positive constant $C_{5}$ which depends on
$\kappa$ and $\eta$  such that for all $t\leq T_{\eta}$, it holds that
\begin{align}\label{eq5.8}
   &\|u^{h}\|_{\mathcal{L}^{\infty}_{t}(\dot{B}^{-1+\frac{3}{p}}_{p,1})}
   +\kappa\|u^{h}\|_{L^{1}_{t}(\dot{B}^{1+\frac{3}{p}}_{p,1})}
   +\|v\|_{\mathcal{L}^{\infty}_{t}(\dot{B}^{-2+\frac{3}{q}}_{q,1})}
   +\kappa\|v\|_{L^{1}_{t}(\dot{B}^{\frac{3}{q}}_{q,1})}\nonumber\\
    &\leq2\big(\|u_{0}^{h}\|_{\dot{B}^{-1+\frac{3}{p}}_{p,1}}+ \|v_{0}\|_{\dot{B}^{-2+\frac{3}{q}}_{q,1}}\big)
   \times \exp\Big\{\int_{0}^{t}(\lambda_{1}\|u^{3}(\tau)\|_{\dot{B}^{1+\frac{3}{p}}_{p,1}}+\lambda_{2}\|u^{3}(\tau)\|_{\dot{B}^{\frac{3}{p}}_{p,1}}^{2}
   +\lambda_{3}\|w(\tau)\|_{\dot{B}^{\frac{3}{r}}_{r,1}})d\tau\Big\}\nonumber\\
   &\leq2\big(\|u_{0}^{h}\|_{\dot{B}^{-1+\frac{3}{p}}_{p,1}}+ \|v_{0}\|_{\dot{B}^{-2+\frac{3}{q}}_{q,1}}\big)\times
   \exp\Big\{C_{5}\big(\|u^{3}_{0}\|_{\dot{B}^{-1+\frac{3}{p}}_{p,1}}^{2}+(\|w_{0}\|_{\dot{B}^{-2+\frac{3}{r}}_{r,1}}+1)
   \exp\big\{C_{5}\|u_{0}^{3}\|_{\dot{B}^{-1+\frac{3}{p}}_{p,1}}\big\}+1\big)\Big\}.
\end{align}
Finally we conclude that if we take $C_{0}$ large enough and $c_{0}$ small enough in \eqref{eq1.9}, then it follows from \eqref{eq5.8} that
\begin{align*}
   \|u^{h}\|_{\mathcal{L}^{\infty}_{t}(\dot{B}^{-1+\frac{3}{p}}_{p,1})}
  +\kappa\|u^{h}\|_{L^{1}_{t}(\dot{B}^{1+\frac{3}{p}}_{p,1})}+\|v\|_{\mathcal{L}^{\infty}_{t}(\dot{B}^{-2+\frac{3}{q}}_{q,1})}
  +\kappa\|v\|_{L^{1}_{t}(\dot{B}^{\frac{3}{q}}_{q,1})}\leq
  \frac{\eta}{2}
\end{align*}
for all $t<T_{\eta}$, which contradicts with the maximality of $T_{\eta}$, thus $T^*=\infty$. We complete the proof of Theorem \ref{th1.2}.
\vskip .3in

\section*{Acknowledgement}
This paper is supported by the National Natural Science Foundation of China (no. 11961030, 12361034) and
 the Natural Science Foundation of Shaanxi Province (no. 2022JM-034).
%\vskip .4in

%\appendix
%\section{}

\end{document}